%% file: kprt2.17.tex
\documentclass{siamltex}
\title{Coarse-graining schemes for stochastic lattice systems with short and long-range interactions}

\author{Markos A. Katsoulakis\thanks{Department of Mathematics,
       University of Massachusetts, Amherst, MA 01003, USA, and Department of Applied Mathematics, 
       University of Crete and  Foundation of Research and Technology-Hellas, Greece, ({\tt markos@math.umass.edu}).}
  \and Petr Plech\'a\v{c}\thanks{Department of Mathematics,
       University of Tennessee, Knoxville, TN 37996, USA, and Joint Institute for Computational Sciences, Oak Ridge National Laboratory,
       ({\tt plechac@math.utk.edu}).}
  \and Luc Rey-Bellet\thanks{Department of Mathematics,
       University of Massachusetts, Amherst, MA 01003, USA, ({\tt lr7q@math.umass.edu}).}
  \and Dimitrios K. Tsagkarogiannis\thanks{
       Universit\'a di Roma, Tor Vergata, Rome, Italy, ({\tt tsagkaro@mat.uniroma2.it}).}
}

\usepackage{amsfonts,amssymb}        %
\usepackage{amsmath}
\usepackage{epsfig}

\input{myppmacro3.tex}
\def\card{\mathrm{card}}
\def\supp{\mathrm{supp}}

\newtheorem{condition}{Condition}[section]

\def\diam{\mbox{diam\,}}
\def\supp{\mbox{supp\,}}
\def\card{\mbox{card\,}}

\def\CGPARQ{Q}
\def\CGPARq{q}

\begin{document}

\maketitle

\begin{abstract}
We develop coarse-graining schemes for stochastic many-particle microscopic models 
with competing short- and long-range interactions on a $d$-dimensional lattice.   
We focus on the coarse-graining of equilibrium Gibbs states and 
using cluster expansions we analyze the corresponding renormalization group map.
We quantify the approximation properties of the coarse-grained terms
arising from different types of interactions and present a hierarchy of
correction terms.  We derive semi-analytical numerical schemes 
that are accompanied with a posteriori  error estimates for coarse-grained 
lattice systems with  short and long-range interactions.
\end{abstract}



\pagestyle{myheadings}
\thispagestyle{plain}
\markboth{M.A. Katsoulakis, P. Plech\'a\v{c}, L. Rey-Bellet, D. K. Tsagkarogiannis}%
	 {Coarse-graining schemes for short and long-range interactions}

\section{Introduction}\label{intro}
Many-particle microscopic  systems with combined short and long-range interactions are  ubiquitous 
in a variety of physical and biochemical systems, 
\cite{andelman}. They exhibit rich mesoscopic and macroscopic morphologies due to competition 
of attractive and repulsive interaction potentials. For example,  
mesoscale pattern formation via self-assembly arises in  heteroepitaxy, \cite{plass},  
other  notable examples include  polymeric systems, \cite{doi},  and
micromagnetic materials, \cite{Hadji}. 
Simulations of such systems rely on molecular methods such as kinetic Monte Carlo (kMC) 
or Molecular Dynamics (MD). However, the presence of long-range interactions severely limits the 
spatio-temporal scales that can be simulated by such direct  computational methods.

On the other hand, an important class of computational tools used for accelerating   
microscopic molecular   simulations  is  the method of {\em coarse-graining}. 
By lumping together degrees of freedom into coarse-grained variables interacting with new, effective potentials
the complexity of the molecular system is reduced, thus yielding accelerated simulation methods 
capable of reaching  mesoscopic length scales.  
Such methods  have been developed for the study and simulation of crystal growth, surface
processes and  polymers,  e.g.,
\cite{KMV, vagelis, kremerplathe,  briels, dionrev07,  KPR}, while there is an extensive literature in soft matter and complex 
fluids, e.g.,  \cite{voth, g219, g97, g98}.
Existing approaches can
give unprecedented speed-up to molecular simulations and can work well
in certain parameter regimes, for instance,  at high temperatures or low densities of the systems. 
On the other hand important macroscopic properties may not be  captured properly in many parameter
regimes, e.g.,
the melt structures of polymers, \cite{kremerplathe};
or the crystallization of complex fluids, \cite{PK06}.  
Motivated in part by such observations
we formulated and analyzed, from a numerical analysis and statistical mechanics perspective,
coarse-grained  variable selection and  error quantification of coarse-grained   approximations 
focusing on stochastic lattice systems with long-range interactions, \cite{KPS, KPRT, AKPR, KPR}. 
We have shown that the  ensuing schemes, 
known as coarse-grained Monte Carlo (CGMC) methods, perform remarkably  well even though traditional 
Monte Carlo methods experience a serious slow-down.
In this paper we focus on lattice systems with both short and long-range interactions. 
Short-range interactions introduce strong correlations between coarse-grained variables 
and a radically different approach needs to be 
employed in order to carry out a systematic and accurate coarse-graining of such systems.

The coarse-graining of microscopic systems is essentially a problem in approximation theory  and 
numerical analysis.   However, the presence of {\em stochastic fluctuations} on 
one hand, and the {\em extensive} nature of the models 
(the presence of extensive quantities that scale as $\BIGO(N)$ with the size of system $N$)
on the other create a new set of challenges.  Before we proceed with the main results of this paper
we  discuss all these issues in a general setting that applies to both on-lattice  and off-lattice 
systems and present the mathematical and numerical framework of coarse-graining
for equilibrium many-body systems.

We denote by $\sigma$  microscopic states of a many-particle system and 
by $\SIGMA$ the set of all microscopic states (i.e., the  configuration space). 
The energy of a configuration is given by the Hamiltonian
$H_N(\sigma)$
where $N$ denotes the size of the microscopic system.
An example studied in this paper is the   $d$-dimensional Ising-type model defined on a lattice 
with $N=n^d$ lattice points, and suitable boundary conditions, e.g., periodic. 
For both on-lattice or off-lattice particle  systems the finite-volume equilibrium states of the system are 
given by the canonical Gibbs measure at the inverse temperature $\beta$, describing the most 
probable configurations
\begin{equation}\label{microGibbs}
    \mu_{N,\beta}(d\sigma)=\frac{1}{Z_N}\EXP{-\beta H_N(\sigma)}
			    P_N(d\sigma)\COMMA
\end{equation}
where the normalizing factor $Z_N=\int e^{-\beta H_N} P_N$,  the partition function, ensures 
that \eqref{microGibbs} is a probability measure,   and $P_N(d\sigma)$ denotes the prior distribution 
on $\SIGMA$. The prior distribution is typically a product measure (see for instance \eqref{microGibbs2}) 
which describes non-interacting particle, or equivalently describes the system at  infinite temperature $\beta=0$. 
At the $\beta=0$ limit the particle interactions included in $H_N$ are unimportant
and thermal fluctuations, i.e., disorder, associated with the product structure of the prior, dominates 
the system.  By contrast at the zero temperature limit, $\beta\to\infty$, interactions dominate and 
thermal fluctuations are unimportant;  in this case \eqref{microGibbs}  concentrates on the minimizers, 
also known as the ``ground states'', of the Hamiltonian $H_N$ over all configurations $\sigma$. 
Finite temperatures, $0<\beta <\infty$, describe intermediate states to these two extreme regimes, 
including  possibly phase transitions, i.e., regimes when as parameters, such as the temperature, 
change, the system exhibits an abrupt transition from a disordered to an ordered state and vice versa, 
or between different ordered phases. 

The objective of (equilibrium) computational statistical mechanics is the simulation of averages over 
Gibbs states, \eqref{microGibbs} of observable quantities $f(\sigma)$
\begin{equation}\label{averageMC}
\EXPECT_{\mu_{N\beta}}[f]\,=\, \int f(\sigma) {\mu}_{N\beta}(d\sigma)\PERIOD
\end{equation}
Due to the exceedingly high dimension of the integration, even for moderate values of the system size $N$,  
e.g., $|\SIGMA|=2^N$ for the standard Ising model, such averaged observables are typically calculated by
Markov Chain Monte Carlo (MCMC) methods, \cite{binder}. 
Nonetheless, mesoscale morphologies, e.g., traveling waves and patterns, are beyond the reach 
of conventional Monte Carlo methods.  For this reason  coarse-graining 
methods have been developed in order to speed up  molecular simulations.

We briefly discuss the mathematical formulation and numerical analysis challenges 
arising in  coarse-graining of   an equilibrium system described by \eqref{microGibbs}.
We rewrite the microscopic configuration $\sigma$ in terms of coarse variables $\eta$ and  
corresponding  fine variables   
$\xi$ so that $\sigma=(\eta, \xi)$. 
We denote the configuration space at the coarse level by $\SIGMAC$
and we denote by $\COP$ the coarse-graining map 
$\COP :\SIGMA \rightarrow \SIGMAC\COMMA\; \COP\sigma=\eta \in \SIGMAC$. 
The coarse-grained system size is denoted by $M$, while the microscopic system size is $N=\CGPARQ M$,
where we  refer to $Q$ as the level of coarse-graining, and $\CGPARQ =1$ corresponds to no coarse-graining.

At the coarse-grained level one is interested in observables  $f(\eta)$ which depend 
only  on the coarse variable $\eta$ and a coarse-grained statistical description of the equilibrium 
properties of the system should be given by a probability measure $\BARIT{\mu}_{M,\beta}(d\eta)$ on 
$\SIGMAC$  such that  the average (the expected value) of such observable is same in the
coarse-grained as well as fully resolved systems. This motivates the following definition. %
\begin{definition} The exact  coarse-grained Gibbs measure $\BARIT{\mu}_{M,\beta}$ is
{\em defined by} 
\begin{equation}\label{ecm00}
\BARIT{\mu}_{M,\beta}( A)  \,\equiv \, \mu_{N,\beta}( \COP^{-1}(A)) \,,
\end{equation}
for any (measurable) set $A\subset \SIGMAC$ or, equivalently,
\begin{equation}\label{ecm0}
\int f(\eta) \,\BARIT{\mu}_{M,\beta}(d\eta) \,=\, \int f(\COP(\sigma)) \,\mu_{N,\beta}(d\sigma) \,.
\end{equation}
for all (bounded) $f: \SIGMAC \to \R$.  
\end{definition}

Slightly abusing notation we will write $\BARIT{\mu}_{M,\beta} \,\equiv \, \mu_{N,\beta} \circ \COP^{-1}$
in the sequel.   In order to write the measure $\BARIT{\mu}_{M,\beta}$ in a more convenient form
we first compute the exact coarse-graining of the prior distribution $P_N(d\sigma)$ on
$\SIGMA$
$$ 
\BARIT{P}_M( d{\eta}) = P_N \circ \COP^{-1} \PERIOD
$$    
The  conditional  prior probability  $P_N(d\sigma\SEP{\eta})$ of having 
a microscopic configuration $\sigma$  given a coarse configuration $\eta$ will  play a crucial  role in the sequel.
Recall that for a function $g(\sigma)$  the conditional  expectation is given by 
\begin{eqnarray}
  \EXPECT[  g \SEP{\eta}] \,&=&\, \int g(\sigma) \,P_N(d\sigma\SEP{\eta})\PERIOD 
\end{eqnarray}
We now write the coarse-grained Gibbs measure $\BARIT{\mu}_{M,\beta}$ using a coarse-grained 
Hamiltonian $\BARH_M (\eta)$.  
\begin{definition} The exact coarse-grained Hamiltonian $\BARH_M(\eta)$ is given by
\begin{eqnarray}\label{rg}
    \EXP{-\beta \BARH_M (\eta)} & \,= \,&\EXPECT[\EXP{-\beta H_N}\SEP{\eta}]\PERIOD
\end{eqnarray}
\end{definition}
This procedure is known as a {\em renormalization group map},  \cite{Kad, Golden}. 
Note that the partition functions for $H_N$ and $\BARH_M$ coincide since
\[
Z_N = \int e^{- \beta H_N} P_N(d\sigma) = \int  \int e^{-\beta H_N} P_N(d\sigma\SEP\eta) \BARIT{P}_M(d\eta)
= \int e^{- \beta \BARH_M} \BARIT{P}_M(d\eta) \equiv \BARIT{Z}_M\PERIOD
\]
Hence for any function $f(\eta)$ we have 
\begin{eqnarray*}
\int f(\eta)\mu_{N,\beta}(d\sigma) \,&=&\, \int f(\eta) \frac{1}{Z_N} e^{- \beta H_N} P_N(d\sigma) 
\,=\, \int f(\eta)\frac{1}{Z_N} \int e^{-\beta H_N} P_N(d\sigma\SEP\eta) \BARIT{P}_M(d\eta) \\ \,&=&\, 
\int f(\eta) \frac{1}{\BARIT{Z}_M} \EXP{-\beta\BARH_M({\eta})} \BARIT{P}_M( d{\eta}) \COMMA
\end{eqnarray*}
and thus the coarse-grained measure  $\BARIT{\mu}_{M,\beta}(d{\eta})$ in \eqref{ecm00} is given by 
\begin{equation}\label{cg_gibbs}
    \BARIT{\mu}_{M,\beta}(d{\eta})=\frac{1}{\BARIT{Z}_M}
    \EXP{-\beta\BARH_M({\eta})} \BARIT{P}_M( d{\eta})\PERIOD
\end{equation}
Although typically $\BARIT{P}_M( d{\eta})$ is  easy to calculate, see e.g., \eqref{microprior}, 
the exact computation of the coarse-grained 
Hamiltonian $\BARH_M(\eta)$ given by \eqref{cg_gibbs} is, in general, an impossible task
even for moderately small values of $N$.

In this paper we restrict our attention to lattice systems, and our  main result is the development 
of a general strategy   to construct explicit  numerical approximations of  the exact 
coarse-grained Hamiltonian $\BARH_M(\eta)$ in the physically important case of combined and 
competing short and long range interactions. Essentially we construct  an approximate 
coarse-grained  energy landscape for the original complex microscopic 
lattice system in Section~\ref{microCG}. 
We show that there is an expansion of   $\BARH_M(\eta)$ into a convergent series 
\begin{equation}\label{series}
\BARH_M(\eta)=\BARH_M^{(0)}(\eta)+\BARH_M^{(1)}(\eta)+\BARH_M^{(2)}(\eta)+\mbox{error}
\end{equation}
by constructing  a suitable  first approximation $\BARH_M^{(0)}(\eta)$ and identifying small 
parameters to control the higher-order terms in the expansion. Truncations including a first few terms in  \eqref{series}
correspond to coarse-graining schemes of increasing accuracy.
In order to obtain this expansion  we rewrite  \eqref{rg} as
\begin{equation}\label{rg2}
\BARH_M(\eta)=\BARH_M^{(0)}(\eta)
-\frac{1}{\beta}\log \EXPECT[\EXP{-\beta (H_N-\BARH_M^{(0)}(\eta))}\SEP{\eta}]\PERIOD
\end{equation}
We need to show that the logarithm can be expanded into a convergent series, uniformly in $N$, yielding eventually an expression 
of the type \eqref{series}. However, two interrelated  difficulties emerge immediately: 
(a) the stochasticity of the system in the finite temperature case yields  the nonlinear expression in \eqref{rg2} 
which in turn will need to be expanded into a series; (b) the extensive nature of the microscopic system, i.e.,
typically the Hamiltonian scales as $H_N=\BIGO(N)$, does not allow the expansion of the logarithm 
and exponential  functions into the Taylor series.

For these reasons, one of the principal mathematical tools we employ
is the  {\em cluster expansion method}, see  \cite{Simon}   for an overview and references. 
As we shall see in the course of this paper  cluster expansions  will allow us to identify 
uncorrelated components in the expected value
$\EXPECT[\EXP{-\beta (H_N-\BARH_M^{(0)}(\eta))}\SEP{\eta}]$,
which in turn will permit us to  factorize it,
and subsequently expand the logarithm in \eqref{rg2} in order to obtain the series \eqref{series}. 
The coarse-graining of systems with purely long-range
interactions was extensively studied using cluster expansions  in \cite{KPRT, AKPR, KPR}. Here we are
broadly following and extending this  approach.  However, the presence of  both short and 
long-range interactions presents new difficulties and requires new methods based on the ideas developed in \cite{op,bco}.  
Short-range interactions induce sub-grid scale correlations between coarse variables, and need to be explicitly 
included in the initial approximation $\BARH_M^{(0)}(\eta)$. To account for these effects we introduce a {\em multi-scale decomposition}  
of the Gibbs state \eqref{microGibbs}
into fine and coarse variables, which in turn  allows us to describe, in a explicit manner, the communication between scales  
for both short and long-range interactions. Furthermore, the multi-scale decomposition of \eqref{microGibbs}
can  also allow us  to  reverse the procedure of coarse-graining in a mathematically systematic manner, i.e.,
{\em reconstruct} spatially localized ``atomistic'' properties,  directly from coarse-grained  simulations. 
We note that this   issue  arises extensively  in the polymer science literature, \cite{tsop, mp}.

An important outcome of the cluster expansion analysis for the approximation 
of \eqref{series}  is the {\em semi-analytical splitting scheme} for the coarse-graining 
of lattice systems with short and long-range interactions. 
Presumably similar strategies could be applied for off-lattice systems such as  the  
coarse-graining of polymers.  The schemes proposed here can be split, within a controllable 
approximation error, into a long and a short-range calculation, see
\eqref{split}.
The long-range part, which is
computationally expensive  for conventional Monte Carlo methods,
can be cheaply  simulated  using  the analytical formula given in
\eqref{l01} in the spirit of our previous work \cite{KPRT}. In this case the saving comes from 
reducing the degrees of freedom  by $\CGPARQ=N/M$ and compressing the range of interactions.
For the short-range interactions  we use the semi-analytical formulas \eqref{semi0} 
which involve precomputing coarse-grained interactions 
with Monte Carlo simulation. However, the simulation is done for  a single  subdomain of  three adjacent coarse cells. 
The error estimates in Theorem~\ref{maintheorem}   also suggest  an improved  decomposition to short and long-range interactions. 
Indeed,  they imply splitting  and rearrangement of  the overall  combined short and long-range potential
into  a new short-range component   that includes  possible singularities originally  in the long-range interaction, e.g.,
the non-smooth part in a Lennard-Jones potential, and  a locally integrable (or smooth)  long-range decaying component. 

In contrast to the splitting approach developed here that allows us to analytically calculate the long range effective Hamiltonian 
(\ref{l02}) in (\ref{split})  and in parallel carry out the semi-analytical   step for (\ref{semi0}), 
existing methods, e.g., (\cite{doi, kremerplathe}), employ semi-analytical  computations involving both short, as well as costly long-range interactions. 
Thus,  multi-body terms, which are believed to be   important at lower temperatures, \cite{doi}, have to be  disregarded.  
A notable result of our  error analysis is the {\it quantification} of the role of multi-body terms in coarse-graining schemes, and the relative 
ease to implement them using the aforementioned splitting schemes.
In Section~\ref{semi}, we further quantify the regimes where such
multi-body  terms are necessary in the context of a specific example.
In \cite{AKPR} the necessity to  include multi-body terms in the effective coarse-grained Hamiltonian  was first discussed in a numerical analysis context 
for systems with  singular (at the origin)  long-range interactions. 

Cluster expansions such as \eqref{series} can also be used  for constructing {\em a posteriori error estimates} for coarse-graining problems, 
based on the rather elementary  observation that
higher-order terms in \eqref{expansion} can be viewed as errors that depend only on the coarse variables $\eta$.
In  \cite{KPRT2} we already employed this type of estimates  for stochastic lattice systems
with long-range interactions in order to construct adaptive coarse-graining schemes.  These tools  operated as an ``on-the-fly''
coarsening/refinement method that recovers accurately phase-diagrams. The estimates allowed us to change adaptively the coarse-graining  
level within the coarse-graining  hierarchy once suitably large or small errors were detected, and thus  to speed up the calculations of phase diagrams. 
Adaptive simulations for molecular systems have been also recently proposed in \cite{clementi}, although they are not based on an  a posteriori error  analysis perspective. 
Finally, the cluster expansions necessary for the rigorous derivation and error estimates of the  schemes developed here rely on the smallness of 
a suitable parameter introduced in Theorem~\ref{maintheorem}, see \eqref{fs-fl}. 
In Section~\ref{semi}, we construct an a posteriori bound for this quantity that can allow 
us to track the validity of the cluster expansion 
for a given resolution in the course of a simulation. 
This approach is, at an abstract level, similar to conditional a posteriori estimates proposed earlier in the numerical analysis 
of geometric partial differential equations, \cite{fierro, nochetto}.

Further challenges for systems with short and long-range interactions  not discussed here include: 
error estimates for observables/quantities  of interest, the development  of coarse-grained dynamics 
from microscopics, phase transitions and estimation of 
physical parameters, such as  critical temperatures. Work related to these  directions for systems with 
long-range interactions have been carried out in \cite{KPS}, \cite{presutti} and  \cite{BZ2}.

The paper is organized as follows. In Section~\ref{microCG} we present  the microscopic Ising-type models with short and long-range interactions
and introduce  the coarse-graining maps and the resulting coarse-grained configuration spaces.
In Section~\ref{motivation} we discuss
our general strategy for the analysis of systems with short and long-range interactions and  
present our main results. 
In Section~\ref{semi} we discuss semi-analytical coarse-graining schemes and their applications to specific examples.  
Section~\ref{polymer} is devoted to the construction of the cluster expansion and to the proof of
convergence of our schemes.

\smallskip\noindent
{\bf Acknowledgments:}
The research of M.A.K. was supported by the National Science Foundation through the grants and NSF-DMS-0715125
and the CDI -Type II award NSF-CMMI-0835673, the U.S. Department of Energy through the grant DE-SC0002339,
and the European Commission   Marie-Curie grant FP6-517911.
The research of P.P. was partially supported by the National Science Foundation under the grant
NSF-DMS-0813893 and by the Office of Advanced Scientific Computing Research,
U.S. Department of Energy under DE-SC0001340; the work was partly done at the Oak Ridge National Laboratory, 
which is managed by UT-Battelle, LLC under Contract No. DE-AC05-00OR22725. 
The research of L. R.-B. was partially supported by the grant NSF-DMS-06058. 
The research of D. K. T. was partially supported by the Marie-Curie grant PIEF-GA-2008-220385. 

\section{Microscopic lattice models  and coarse-graining}\label{microCG}
We consider an Ising-type model on the $d$-dimensional  square lattice 
$\LATT:=\{ x=(x_1, \cdots, x_d) \in \Z^d \,;\, 0\le x_i \le n-1\}$ with 
$N=n^d$ lattice points.   For simplicity we assume periodic boundary 
conditions throughout this
paper although other boundary conditions can be accommodated. 
At each lattice site $x$ there is a spin $\sigma(x)$ taking values in
$\Sigma=\{+1,-1\}$. A spin configuration $\sigma=\{\sigma(x)\}_{x\in\LATT}$
on the lattice $\LATT$ is an element of the configuration space
$\SIGMA:=\Sigma^{\LATT}$.  For any subset $X\subset\LATT$ we denote 
$\sigma_{X}=\{\sigma(x)\}_{x\in X}\in\Sigma^{X}$ the restriction of 
the spin configuration to $X$. Similarly, for a function $f :\SIGMA\to\R$ we
denote $f_X$ the restriction of $f$ to $\Sigma^X$.
The energy of a configuration $\sigma$ is given by the Hamiltonian
\begin{equation}\label{microHamiltonian}
H_N(\sigma) = H_N^{s}(\sigma) + H_N^{l}(\sigma) \,, 
\end{equation}
which consists of a short-range part $H_N^{s}$ and a long range part $H_N^{l}$. 
For the short-range part we have
\[
H_N^{s}(\sigma) \,=\, \sum_{X\subset \LATT}U_X(\sigma)\COMMA
\]
where the {\it short-range} potential $U=\{ U_X,\,X\subset \Z^d\}$, with $U_X: \Sigma^X\to\R$, 
is translation invariant (i.e., $U_{X+y}=U_X$ for all $X \subset \Z^d$ and all $y \in Z^d$) 
and has the finite range $S$ (i.e., $U_X=0$ whenever $\diam(X)> S$). 
We define the norm $\|U\|\equiv \sum_{X\supset\{0\}\SEP\mathrm{diam}(X)\le S} \|U_X\|_\infty$ where the norm
$\|\cdot\|_\infty$ is the standard sup-norm on the space of continuous functions.
A typical case is the  nearest-neighbor Ising model
\[
H^{s}_N(\sigma)=K\sum_{\langle x,y\rangle}\sigma(x)\sigma(y)\COMMA
\] 
where by $\langle x,y\rangle$ we denote summation over the nearest neighbors. 
For the {\it long-range part} we assume the form 
\[
H_N^{l}(\sigma) \,=\, -\frac{1}{2}\sum_{x\in\LATT}\sum_{y\neq x}J(x-y)\sigma(x)\sigma(y) \COMMA
\]
where the two-body potential $J$ has the form 
\[
J(x-y) = \frac{1}{\RANGE^d}\JNOT\left(\frac{1}{\RANGE}|x-y|\right)\COMMA
            \SPACE \label{defJV1}\\
\]
for some $V \in C^1 ([0,\infty))$.  The factor $1/\RANGE^d$ in
\eqref{defJV1} is a normalization which ensures that the strength of
the potential $J$ is essentially independent of $L$, i.e., $\sum_{x \ne
0} |J(x)| \simeq \int |V(r)| dr$.  For example, if we choose $V$ such
that $V(r)=0$ for $r>1$ then a spin at the site $x$ interacts with its
neighbors which are at most $\RANGE$ lattice points away from $x$ and
in this case $L$ is the range of the interaction $J$.  
It is convenient to think of 
$L$ as a parameter in our model and  more precise assumptions on the interactions will 
be specified later on. 

The finite-volume equilibrium states of the system are given by the
canonical Gibbs measure \eqref{microGibbs}
and $P_N(d\sigma)$, the prior distribution on $\SIGMA$, is a product measure
\begin{equation}\label{microGibbs2}
 P_N(d\sigma)=\prod_{x\in\LATT} P_x(d\sigma(x))\PERIOD
\end{equation}
A typical choice is $P_x(\sigma(x)=+1)=\frac{1}{2}$ and $P_x(\sigma(x)=-1)=\frac{1}{2}$, 
i.e., independent Bernoulli random variables at each site $x\in\LATT$.
For the sake of simplicity we consider Ising-type spin systems, but  
the techniques and ideas in this paper apply also to  Potts and  
Heisenberg models or, more generally, to models where the ``spin'' 
variable  takes values in a compact space. 

\subsection{Coarse-graining} 
In order to coarse-grain our system we divide the lattice $\LATT$ into coarse
cells and define coarse variables by averaging spin values over the coarse 
cells.  We partition the lattice $\LATT$ into $M=m^d$ disjoint cubic coarse cells, 
each cell containing $\CGPARQ=q^d$ microscopic lattice points so that $N = n^d=(m
q)^d = M\CGPARQ$. The coarse-grained (real-space) hierarchy  can be build in a anisotropic
way, by replacing $n$, $m$, $q$ with multi-indexes. For example, different levels of
coarse-graining in individual coordinate directions will be given by $q=(q_1,\dots,q_d)$ and
the power $q^d$ would be interpreted as $q_1q_2\dots q_d$. We refrain from an unnecessary 
generality and assume that the coarse-graining is isotropic, $q_1 = \dots = q_d = q$.
We define a coarse lattice $\LATTC = \{ k=(k_1,\cdots,k_d) \in \Z^d
\,;\, 0 \le k_i < m-1 \}$ and we set $\LATT = \cup_{k \in \LATTC}
\CUBE_k$ where $\CUBE_k =\{ x \in \LATTC \,;\, k_i q \le x_i <
(k_i+1)q\}$.  Whenever convenient we will identify the coarse cell
$\CUBE_K$ in the microscopic lattice $\LATT$ with the point $k$ of the
coarse lattice $\LATTC$.
For any configuration $\sigma^k \equiv \sigma_{C_k}$ on the coarse cell $\CUBE_k$ we assign a new spin value  
$$
\eta(k)=\sum_{x\in \CUBE_k}\sigma(x)
$$ 
which takes values in $\bar{\Sigma}=\{ -\CGPARQ, -\CGPARQ+2, \ldots, \CGPARQ\}$. 
We denote the configuration space at the coarse level by $\SIGMAC \equiv \bar{\Sigma}^{\LATTC}$
and we denote by $\COP$ the coarse-graining map 
$$
 \COP :\SIGMA \rightarrow \SIGMAC\COMMA\;\;\;\; 
       \sigma=\{\sigma(x)\}_{ x\in \LATT} \mapsto \eta = \{ \eta(k)\}_{  k \in \LATTC }
$$ 
which assigns a configuration $\eta$ on the coarse lattice $\LATTC$ given a
configuration $\sigma$ on the microscopic lattice $\LATT$. 

The exact   coarse-grained Gibbs measure is defined in \eqref{ecm00} for arbitrary Gibbs states having the form \eqref{cg_gibbs}.
Since $\eta(k)$ depends only on the spins $\sigma(x)$, with
$x \in \CUBE_k$, the coarse-grained measure ${\bar P}_M$ is a product measure
\begin{equation}\label{microprior}
\BARIT{P}_M( d{\eta}) = P_N \circ \COP^{-1} = \prod_{k \in \LATTC} \bar{P}_k (d \eta(k))\PERIOD
\end{equation}
For example if $P_x$ is a Bernoulli distribution then 
$P_{k}(\eta(k)) = \binom{\CGPARQ}{\tfrac{\eta(k)+ \CGPARQ}{2}}\left(\tfrac{1}{2}\right)^{\CGPARQ}$.
 Similarly, we define the  conditional probability measure $P_N(d\sigma\SEP{\eta})$ of having 
a microscopic configuration $\sigma$ on $\LATT$ given a coarse configuration $\eta$ on 
$\LATTC$. This measure  plays a crucial  role in the sequel since it factorizes over the coarse cells
\begin{equation}\label{cond_meas}
   P_N(d\sigma\SEP{\eta})= \prod_{k \in \LATTC} P_k(d\sigma^{k} \SEP \eta(k))\COMMA
\end{equation}
where $P_k(d\sigma^{k} \SEP \eta(k))$ is the conditional probability of a microscopic configuration 
$\sigma^k$ on $\CUBE_K$ given a coarse configuration $\eta(k)$. 

\section{Approximation strategies for $\BARH_M(\eta)$}\label{motivation}
In this section we present a general strategy for constructing
approximations of  the exact coarse-grained Hamiltonian $\BARH_M(\eta)$ in \eqref{cg_gibbs}.  
We show how to expand $\BARH_M(\eta)$ into a convergent series \eqref{series}
by choosing a suitable first approximation $\BARH_M^{(0)}(\eta)$ and identifying small 
parameters to control the higher-order terms in the expansions.  The basic idea
is to use the first approximation $\BARH_M^{(0)}(\eta)$ in order to rewrite \eqref{rg} as
\eqref{rg2}.
We show that the logarithm can be expanded into a convergent series, uniformly in $N$, 
using suitable cluster expansion techniques.   
We discuss in detail the case $d=1$ in order to illustrate general
ideas in the case where calculations and formulas are relatively simple.  The general
$d$-dimensional case
is discussed in detail in Section~\ref{polymer}.

We recall that the Hamiltonian $H_N(\sigma)=H_N^{l}(\sigma)+H_N^{s}(\sigma)$ consists
of a short-range part $H_N^{s}(\sigma)$ with the range $S$ and a
long-range part $H_N^{l}(\sigma)$ whose range is  $L$. 
We choose the coarse-graining level $\CGPARq$ such that
$$
S \,<\,  \CGPARq \, < \,L  \PERIOD
$$ 
There are two small parameters associated with the range of the interactions 
$$
\epsilon_s \propto  \frac{S}{q}\,, \quad\quad\mbox{and}\quad\quad \epsilon_l \propto  \frac{q}{L}\PERIOD
$$ 
The first approximation is of the form 
\begin{equation}\label{split0}
\BARH_M^{(0)} = \BARH_M^{l,(0)} + \BARH_M^{s,(0)}\COMMA
\end{equation}
and two distinct separate procedures are used to define 
the short-range coarse-grained  approximation $\BARH_M^{s,(0)}$, as well as its   long-range counterpart 
 $\BARH_M^{l,(0)}$. Due to the nonlinear nature of the map induced by \eqref{rg2} it is not obvious that \eqref{split0} will 
be a valid approximation, except possibly at high temperatures, when $\beta <\!< 1$. This fact will be established for a wide range of parameters
in  the error analysis of Theorem~\ref{maintheorem}, and in the discussion in Section~\ref{semi}, provided  a suitable choice is made for  
$\BARH_M^{s,(0)}$ and  $\BARH_M^{l,(0)}$.

\subsection{Coarse-graining of the long-range interactions}  We briefly recall the coarse-graining 
strategy of \cite{KPRT} for the long-range interactions. Since the range of the interaction, $L$, is 
larger than the range of coarse-graining $\CGPARQ$ a natural first approximation for the long-range part is to 
average the interaction $J(x-y)$ over coarse cells.  Thus we define  
\begin{equation}\label{l01}
\BARH_M^{l,(0)}(\eta) \,\equiv E[ H_N^{l}\SEP\eta]\COMMA
\end{equation}
and an easy computation gives  
\begin{equation} \label{l02}
\BARH_M^{l,(0)}(\eta) \,=\, \, - \frac{1}{2} \sum_{k \in \LATTC} \sum_{l\not= k} \bar{J}(k,l) \eta(k) \eta(l) -
                                 \frac{1}{2} \sum_{k \in \LATTC} \bar{J}(k,k) ( \eta(k)^2 - \CGPARQ)  \COMMA
\end{equation}
where
\[
\bar J(k, l)=\frac{1}{\CGPARQ^2} \sum_{x\in \CUBE_k}\sum_{y\in \CUBE_l} J(x-y) \,, \quad 
\bar J(k, k)=\frac{1}{\CGPARQ(\CGPARQ-1)} \sum_{x, y \in \CUBE_k}\sum_{y\neq x } J(x-y) \PERIOD
\]
A simple error estimate (see \cite{KPRT,AKPR} for details in various cases) gives 
\[
H_N^{l}(\sigma)=\BARH_M^{l,(0)}( \COP(\sigma)) + e_L \quad {\rm with } \quad 
e_L=N \BIGO(\frac{\CGPARq}{L} \|\nabla V\|_\infty)\PERIOD
\]
Using this definition of $\BARH_M^{l,(0)}$ we obtain
\begin{eqnarray}
\label{l04}
\EXP{-\beta H_N^{l} (\sigma)  } P_N(d\sigma \SEP \eta) &=&  \EXP{-\beta \BARH_M^{l,(0)}(\eta) }  
                                                            \EXP{-\beta \left[ H_N^{l}(\sigma) - \BARH_M^{l,(0)}(\eta) \right] } 
          P_N(d\sigma\SEP \eta)\COMMA \label{l05}  \\
&=&  \EXP{-\beta \BARH_M^{l,(0)}(\eta) }  \prod_{j, k  \in \LATTC} \left(1 + f^l_{jk} \right)  P_N(d\sigma\SEP \eta)\COMMA \label{pol1l}
\end{eqnarray} 
where  
\begin{equation}\label{flong}
f^l_{jk} \equiv  \EXP{ \frac{\beta}{2} \sum_{x\in \CUBE_j} \sum_{y\in \CUBE_k,y\neq x}
   (J(x-y)-\bar{J}(k,l))\sigma(x)\sigma(y)(2-\delta_{jk})}  -1   \PERIOD
\end{equation}
Due to the fact that $P_N(d\sigma\SEP\eta)$ has a product structure one can rewrite \eqref{pol1l} 
as a cluster expansion, \cite{KPRT} (see also Section~\ref{polymer}), 
as in \eqref{series}. The key element in that cluster expansion is the ``smallness'' of the quantity
\begin{equation}\label{errJb}
  |J(x-y) - \bar J(k,l)| \leq 2\frac{\CGPARq}{L^{d+1}}
          \sup_{\genfrac{}{}{0cm}{2}{x'\in \CUBE_k,}{y'\in \CUBE_l}}|\nabla V(x'-y')|\, ,
\end{equation}
which yields asymptotics
\begin{equation}\label{smallness0}
f^l_{jk} \sim
   \BIGO(\CGPARq^{2d}\frac{\CGPARq}{L^{d+1}}\|\GRADV\|_\infty)\PERIOD
\end{equation}  
The estimate \VIZ{errJb} follows from regularity assumptions on $V$ and the Taylor expansion. 

\subsection{Coarse-graining of short-range interactions}  For the short-range part, using  that $S < q$, 
we write the Hamiltonian as 
\begin{equation}\label{decompH}
H_N^{s}(\sigma)=\sum_{k\in\LATTC} H_k^{s}(\sigma)+
\sum_{k\in\LATTC} W_{k,k+1}(\sigma)\COMMA
\end{equation}
where
\[
H_k^{s}(\sigma)=\sum_{X \subset \CUBE_k} U_X(\sigma) \COMMA \quad  W_{k,k+1}(\sigma)= \sum_{X \cap \CUBE_{k} \not= \emptyset \COMMA\, 
X \cap \CUBE_{k +1} \not= \emptyset} U_X(\sigma) \COMMA
\] 
i.e., $H^{s}_k$ is the energy for the cell $\CUBE_k$ which does not interact with other cells, i.e., under 
the free boundary conditions, and 
$W_{k,k+1}$ is the interaction energy between the cells $\CUBE_k$ and $\CUBE_{k+1}$.
Note the elementary bound 
\begin{equation}\label{w}
\sup_{\sigma}W_{k,k+1}(\sigma) \sim S q^{d-1}\|U\|\PERIOD
\end{equation}
The most naive  coarse-graining, besides of course developing a mean-field-type  approximation, consists 
in regarding the boundary terms $W_{k,k+1}$ as a perturbation.  We have then, formally,   
\begin{eqnarray*}
\EXP{- \beta \BARH_M(\eta)} & \sim & \int \EXP{- \beta \BARH_M^{l,(0)}(\eta)+e_L+e_S}e^{-\sum_{k\in\LATTC} \beta H_{\CUBE_k}^{s}(\sigma)} P_N(d\sigma\SEP\,\eta)\\
& = & \EXP{-\beta\BARH_M^{l,(0)}(\eta)+e_L+e_S}\prod_{k\in\LATTC} e^{-\beta \bar{U}_k^{s,(0)}(\eta_k)}\COMMA
\end{eqnarray*}
where the one-body potential 
$$
\bar{U}^{s,(0)}_k(\eta_k) = -\frac{1}{\beta} \log \int e^{-\beta H_k^{s}(\sigma)} P_k(d\sigma^k | \eta(k)) 
$$ 
is the exact  coarse-grained Hamiltonian for the cell $\CUBE_k$ with free boundary conditions. 
As a result an initial   guess for the zero order approximation could be 
\begin{equation}\label{naive}
\BARH_M^{l,(0)}(\eta)+\sum_k \bar{U}^{s,(0)}_k(\eta_k) \PERIOD
\end{equation}
However,  this approach appears to be rather  simplistic in general  since the
correlations between the cells induced by the short-range potential have
been completely ignored.  While this approximation may be reasonable
at  high temperatures it is not a good starting point for a series
expansion of the Hamiltonian using a cluster expansion.   
Instead we need to adopt a more systematic approach outlined in the next section.

\subsection{Multiscale  decomposition of Gibbs states} 
This approach provides  the common underlying structure of all coarse-graining schemes at equilibrium including lattice and off-lattice models.
It is essentially a  decomposition of the Gibbs state \eqref{microGibbs} into  product measures among different scales selected with suitable properties. 
We outline it for the case of short-range interactions where we rewrite the Gibbs measure \eqref{microGibbs} as
$$
\mu_{N,\beta}(d\sigma) \sim \EXP{-\beta H_N(\sigma)}
                            P_N(d\sigma)=e^{-\beta H_N(\sigma)}P_N(d\sigma\SEP\eta)\BARIT{P}_M(d\eta)\PERIOD
$$
We use the notation $\sim$ meaning up to a normalization constant, i.e., in the equation above we
do not spell out the presence of the constant $Z_N$.
We now seek the following decomposition of the short-range interactions
\begin{equation}\label{general}
\EXP{-\beta H_N^s(\sigma)}P_N(d\sigma\SEP\eta) =R(\eta)A(\sigma)\nu(d\sigma|\eta)\COMMA
\end{equation}
where

\smallskip
\noindent
(a) $R(\eta)$ depends only on the coarse variable $\eta$ and is related to the first coarse-grained approximation $\BARH_M^{(s,0)}(\eta)$ via the formula 
\begin{equation}\label{abstractdecomp}
R(\eta) = \EXP{-\beta   \BARH_M^{s,(0)}(\eta)}\COMMA \quad 
A(\sigma)\nu(d\sigma|\eta)=\EXP{-\beta \left(H_N^s(\sigma) - \BARH_M^{s,(0)}(\eta)\right)}P_N(d\sigma\SEP\eta)\COMMA
\end{equation}

\noindent
(b) $A(\sigma)$ has a form amenable to a cluster expansion, i.e., for $d=1$
\begin{equation}\label{decomp2}
A(\sigma) = \prod_{k\in K} ( 1+ \Phi_k(\sigma)) 
\end{equation}
for some $K\subset \bar{\Lambda}_M$. The function 
$\Phi_k$ is {\em small} and moreover $\Phi_k(\sigma)$ depends on 
the configuration $\sigma$ only locally, up to a fixed finite distance from $\CUBE_k$. 
In the example at hand (for $d=1$) we have $\Phi_k(\sigma)=\Phi_k( \sigma^{k-1}, \sigma^{k+1})$.  

\smallskip
\noindent 
(c) The measure $\nu(d\sigma|\eta)$ has the general form 
\begin{equation}\label{decomp3}
\nu(d\sigma|\eta) = \prod_{k\in\bar{\Lambda}_M}\nu_k(d\sigma|\eta)\COMMA
\end{equation}
where $\nu_k (d\sigma \SEP \eta )$ depends on $\sigma$ and $\eta$ only locally up 
to a fixed finite distance from $\CUBE_k$.  In the example at hand $\nu_k (d\sigma \SEP \eta )$
depends only on the configuration on $\CUBE_{k-1} \cup \CUBE_k \cup \CUBE_{k+1}$. 
Even though the measure $\nu(d\sigma|\eta)$ is not a product measure, the fact that this measure 
has  finite spatial correlation makes it adequate for a cluster expansion, 
see  \eqref{nu} and Section~\ref{polymer}.

Although here we described the multiscale decomposition of the Gibbs measure  for the case of short-range interactions,   
the results on the long-range interactions, discussed earlier, can be reformulated in a similar way.  
In particular, \eqref{l04} and \eqref{l05} can be rewritten as 
\begin{equation}\label{general-long}
\EXP{-\beta H_N^{l}(\sigma)}P_N(d\sigma\SEP\eta) = R(\eta)A(\sigma)\nu(d\sigma|\eta)\COMMA
\end{equation}
where
$R(\eta) = \EXP{-\beta   \BARH_M^{l,(0)}(\eta)}$, $ \nu(d\sigma|\eta)=P_N(d\sigma\SEP\eta)$,
and 
\begin{equation}\label{abstractdecomp-long}
A(\sigma)=\EXP{-\beta\big(H_N^l(\sigma)\beta -\BARH_M^{l,(0)}(\eta)\big)}
         =\prod_{j, k  \in \LATTC} \left(1 + f^l_{jk}\right)\PERIOD
\end{equation}
We recall that in analogy to \eqref{decomp3}, the product structure of $\nu(d\sigma|\eta)=P_N(d\sigma\SEP\eta)$ allows us 
to carry out a cluster expansion for the long-range case, and obtain a convergent series such as \eqref{series}, 
thus yielding an expansion  of the exact coarse-grained Hamiltonian $\BARH_M^{l}$, \cite{KPRT}.
 
We  note that \eqref{general}, used here as a numerical and multiscale analysis tool in order to derive suitable 
approximation schemes for the coarse-grained Hamiltonian,  was  first introduced in \cite{olivieri,op, bco}
for  the purpose of deriving  cluster expansions for lattice systems with short-range interactions
 away from the well-understood high temperature regime.

\subsection{Coarse-graining schemes in one spatial dimension}
We sketch how to obtain a decomposition such as \eqref{general} for $d=1$ and construct suitable $R(\eta)$.
We split the one-dimensional lattice  
into non-communicating components, for instance, even- and odd-indexed cells and write 
\begin{eqnarray}
\EXP{-\beta H_N^s} P_N(d\sigma\SEP\eta) = && \prod_{k:\,\mathrm{odd}}\left[
\EXP{-\beta (W_{k-1,k}+W_{k,k+1})}e^{-\beta H^s_{k}} P_k(d\sigma^k\SEP \eta(k))
\right] \times \nonumber \\
&& \prod_{k:\,\mathrm{even}}\EXP{-\beta H^s_{k}} P_k(d\sigma^k\SEP \eta(k)) \PERIOD \label{sr1}
\end{eqnarray}
In \eqref{sr1} we will normalize the factors for $k$ odd by dividing each factor with the 
suitably defined corresponding partition functions  for the regions  $\CUBE_k$ and $\CUBE_{k-1}\cup\CUBE_k\cup\CUBE_{k+1}$.

\begin{definition} We define the  partition function
with  boundary conditions $\sigma^{k-1}$ and 
$\sigma^{k+1}$, i.e.,
\begin{equation}\label{zboundary}
Z_k(\eta(k); \sigma^{k-1},\sigma^{k+1})\,=\, \int  \EXP{-\beta (W_{k-1,k}+W_{k,k+1})} e^{-\beta H^s_k} P_k(d\sigma^k\SEP \eta(k)) \PERIOD
\end{equation}
In order to decouple even and odd cells we define  the partition function with free boundary conditions on $\CUBE_{k-1}$
and boundary condition $\sigma^{k+1}$ on $\CUBE_{k+1}$, i.e., 
\begin{equation}\label{z0boundary}
Z_k(\eta(k); 0,\sigma^{k+1})\,=\, \int  \EXP{-\beta W_{k,k+1}} \EXP{-\beta H^s_{k}} P_k(d\sigma^k\SEP \eta(k))\COMMA
\end{equation} 
and similarly $Z_k(\eta(k); \sigma^{k-1}, 0)$, as the partition function with free boundary conditions on $\CUBE_{k+1}$ and boundary condition 
$\sigma^{k-1}$ on $\CUBE_{k-1}$.  We also denote by $Z_k(\eta(k); 0,0)$ the partition function for $\CUBE_k$ with free boundary conditions. 
We define  the three-cell partition function with free boundary conditions
\begin{eqnarray}
\label{z00boundary}
&&Z_{k-1,k,k+1}(\eta(k-1),  \eta(k), \eta(k+1); 0,0)=\nonumber\\
\nonumber\\
&&\quad\int \EXP{- \beta \left(H^s_{k-1}+ W_{k,k-1} + H^s_k +  W_{k,k+1} + H^s_{k+1}\right) }\times \nonumber\\
&&\quad\quad
P_{k-1}(d\sigma^{k-1}\SEP \eta(k-1)) P_k(d\sigma^k\SEP \eta(k)) P_{k+1}(d\sigma^{k+1}\SEP \eta(k+1))  \PERIOD
\end{eqnarray}
\end{definition}
The key to the decomposition and eventually to the cluster expansion is the introduction of a ``small term'' analogous to \eqref{smallness0}. 
\begin{definition}
\begin{equation}\label{smallness}
 f^s_{k-1, k+1}(\eta(k); \sigma^{k-1}, \sigma^{k+1})  \,=\,  \frac{Z_k(\eta(k); \sigma^{k-1},\sigma^{k+1})Z_k(\eta(k); 0,0)}
{Z_k(\eta(k); 0,\sigma^{k+1})Z_k(\eta(k); \sigma^{k-1},0)}-1
\end{equation}
\end{definition}
An important element in the cluster expansion in Section~\ref{polymer} is the 
estimation of the terms $f^s_{k-1,k+1}$.
However, a straightforward  
estimate based on \eqref{w} would yield
\begin{equation}\label{asymptotics}
 f^s_{k-1, k+1}(\eta(k); \sigma^{k-1}, \sigma^{k+1})  \,\sim\, \beta S\|U\|\PERIOD
\end{equation}  
We rewrite
\begin{eqnarray}
Z_k(\eta(k);\sigma^{k-1},\sigma^{k+1})&=&
\left( f_{k-1, k+1}(\eta(k); \sigma^{k-1}, \sigma^{k+1})  +1 \right)\times\nonumber \\
&& \quad\quad \frac{Z_k(\eta(k); 0,\sigma^{k+1})Z_k(\eta(k); \sigma^{k-1},0)}{Z_k(\eta(k);0,0)}\PERIOD \label{dec1}
\end{eqnarray}
In \eqref{sr1} we now divide and multiply each factor with $k$ odd by $Z_k(\sigma^{k-1},\sigma^{k+1})$ 
and use the formula \eqref{dec1}. Furthermore,
we multiply each factor with even $k$ by $Z_{k-1,k,k+1}(0,0)$ and obtain 
\begin{eqnarray}
&&\EXP{-\beta H_N^s}P_N(d\sigma\SEP\eta) =  \nonumber \\
&&\underbrace{ \prod_{k:\,\mathrm{odd}} Z_k(0,0)^{-1} \prod_{k:\,\mathrm{even}} Z_{k-1,k,k+1}(0,0)}_{ \displaystyle \equiv R(\eta) }  \,\,
\underbrace{ \prod_{k:\,\mathrm{odd}} (f^s_{k-1,k+1} + 1)}_{\displaystyle \equiv A(\sigma)} \times \label{aa}  \\
&& \underbrace{ 
\prod_{k:\,\mathrm{odd}}  \frac{ e^{-\beta \big( H^s_k + W_{k-1,k} + W_{k,k+1}\big)}}{ Z_k(\sigma^{k-1},\sigma^{k+1})} 
P_k(d \sigma^k\SEP \eta(k))  
\prod_{k:\,\mathrm{even}} \frac{e^{-\beta H^s_k}  Z_{k+1}(\sigma^k , 0) Z_{k-1}(0,\sigma^k)} {Z_{k-1,k,k+1}(0,0)}P_k(\sigma^k\SEP \eta(k))}_{\displaystyle  \equiv \nu(d\sigma|\eta))}\nonumber \\
\label{nu}
\end{eqnarray}
where we have used that 
$$
\prod_{k:\,\mathrm{odd}} Z_k(0,\sigma^{k+1})Z_k(\sigma^{k-1},0) = \prod_{k:\,\mathrm{even}} 
Z_{k+1}(\sigma^k, 0) Z_{k-1}(0, \sigma^k)\PERIOD
$$
It is easy to verify that $\nu(d\sigma\SEP\eta)$ defined in \eqref{nu} 
is a normalized measure and has the form required in condition (c) of the multiscale decomposition of the Gibbs measure. 
The factor $R(\eta)$ defined in \eqref{aa} gives the first order corrections induced by the correlations 
between adjacent cells. Putting together the analysis for short and long-range interactions we obtain  
the main result formulated as a theorem.

\begin{theorem}\label{maintheorem}
Let
\begin{equation}\label{split}
\BARH_M^{(0)}(\eta) \,=\,   \BARH_M^{l, (0)}(\eta)  +    \BARH_M^{s, (0)}(\eta)
\end{equation}
where $\BARH_M^{l, (0)}(\eta)$ is given in \eqref{l01} and \eqref{l02}  and 
\begin{equation}\label{decomp}
\BARH_M^{s, (0)}(\eta) \,=\,  \sum_{k:\,\mathrm{odd}}  \bar{U}_k^{s,(0)}(\eta(k))  + 
 \sum_{k:\,\mathrm{even}}  \bar{U}_{k-1,k,k+1}^{s,(0)}(\eta(k-1), \eta(k), \eta(k+1)) \COMMA
\end{equation}
with the one-body interactions
\begin{equation}\label{part1}
\bar{U}_k^{s,(0)}(\eta(k)) = -\frac{1}{\beta}\log Z_k(\eta(k); 0,0)  \,,
\end{equation}
and the three-body interactions 
\begin{eqnarray}
\bar{U}_{k-1,k,k+1}^{s,(0)}(\eta(k-1), &&\eta(k), \eta(k+1))=  \nonumber
\\
&&-\frac{1}{\beta}\log Z_{k-1,k,k+1}(\eta(k-1), \eta(k), \eta(k+1); 0, 0) \,, \label{part2}
\end{eqnarray}
where $Z_k$ and $Z_{k-1,k,k+1}$ are given by (\ref{z0boundary}) and (\ref{z00boundary}) 
respectively. Then

\noindent
1. we  have the error bound 
$$
|\BARH_M - \BARH_M^{(0)}| \sim N\BIGO\left(\frac{\beta S\|U\|} {\CGPARq}+ \frac{\CGPARq\beta \|\nabla V\|_\infty}{L}\right)\COMMA
$$
for a short-range potential with the range $S<\!< q<\!<L$. 
The loss of information when coarse-graining at the level $\CGPARq$ is quantified by
 the specific relative entropy error 
\begin{equation}\label{relent}
\frac{1}{N}\RELENTR(\BARIT{\mu}^{(0)}_{M,\beta}\SEP\mu_{N,\beta}\circ\COP^{-1}) \,=\, 
\BIGO\left(\frac{\beta S\|U\|} {q}+ \frac{q\beta \|\nabla V\|_\infty}{L}\right)\PERIOD
\end{equation}

\noindent
2. There exist  $\delta_0>0$ and $\delta_1>0$ such that 
if 
\begin{equation}\label{fs-fl}
\sup_k\sup_{\sigma^{k-1}, \sigma^{k+1}, \eta(k)}
|f^s_{k-1, k+1}(\eta(k); \sigma^{k-1}, \sigma^{k+1})| \le \delta_0\, , \quad
\sup_{k, j}\sup_{\sigma^{j}, \sigma^{k}}
|f^l_{jk}(\sigma^{j}, \sigma^{k})| \le \delta_1\COMMA 
\end{equation}
where  $f^s_{k-1, k+1}$ and $f^l_{jk}$ are  given by (\ref{smallness}) and  (\ref{flong}) respectively, then 
$\BARH_M - \BARH_M^{(0)}$ is expanded in a convergent series in the parameter 
$\delta \sim \big(\frac{\beta \|U\|  S} {q} + \frac{q\beta \|\nabla V\|_\infty}{L}\big)$
\begin{equation}\label{expansion}
\bar{H}_M(\eta) =  \bar{H}_M^{(0)}(\eta) + \bar{H}_M^{(1)}(\eta)  + \cdots +
\bar{H}_M^{(p)}(\eta) + M\BIGO(\delta^{p+1})\, .
\end{equation}
\end{theorem}

\medskip
\begin{remark}
{\rm
The error estimate \VIZ{relent} suggests qualitatively an estimate on 
the regimes of validity of the method, and on the
 ``optimal'' level, $q=q_{\mathrm{opt}}$, when we restrict to the regime $S<q<L$, 
where $S$ and $L$ are the respective interaction ranges for short and long-range potentials. 
The corresponding error is then
\begin{equation}\label{optq1}
  q_{\mathrm{opt}} \sim \sqrt{ SL\frac{\|U\|}{\|\nabla V\|_{\infty}}}
    \COMMA \quad 
  \quad
  \frac{1}{N}\RELENTR(\BARIT{\mu}^{(0)}_{M,\beta}\SEP\mu_{N,\beta}\circ\COP^{-1})=
     \BIGO\left(\beta \sqrt{\frac{S}{L} 
        \|U\| \|\nabla V\|_\infty}\right)\, .
\end{equation}
}
\end{remark}

The application of Theorem~\ref{maintheorem}   requires to check the validity of \eqref{fs-fl}.
 Certainly the conditions \eqref{smallness0} and \eqref{asymptotics} 
are satisfied in suitable regimes, see also Section~\ref{polymer} for more details. 
More interestingly, for specific examples these conditions can be verified directly,
we refer to Section~\ref{semi}. In particular, in \eqref{conditional} and \eqref{conditional3}
we even obtain an upper bound that depends only on the coarse observables. 
This allows us to check the conditions \eqref{fs-fl} (dictated by the cluster expansions)  computationally  
in the process of a Monte Carlo simulation involving only the coarse variables $\eta$.

On the other hand, 
in \cite{olivieri,op},  the short-range  condition in \eqref{fs-fl} is taken as an 
assumption.
In one dimension,  this condition  holds up to 
very low temperatures while in dimension $d\ge 2$ this condition can be 
satisfied in the high-temperature regime, see for example the analysis in \cite{bco}
where similar conditions are used for the nearest-neighbor Ising model in the
dimension $d=2$ all the way up to the critical temperature. 

Finally, we note that  a similar strategy to coarse-grained short and long-range interactions can be
used in any dimension, as  we discuss in  Section~\ref{polymer}.  
In the multi-dimensional case we split the domain into boxes of size larger than the range
of the interaction so that the next-to-nearest coarse cells are independent. In one dimension, 
this procedure gives rise to the separation into odd- and even-indexed coarse cells, while     
in higher dimensions it is done in a recursive manner, proceeding one
dimension at a time.  
Then by freezing the configurations on the collection of independent coarse cells
(resulting to the one-body coarse-grained terms)
we create further correlations which couple the remaining cells. 
This fact in one-space dimension yields  the three-body terms, noting 
that possible two-body coarse-grained correlations are contained therein, see also  \eqref{newformula}. 
We also remark that coarse-graining schemes for the nearest-neighbor Ising model, 
involving only two-body interactions were recently proposed in \cite{sinno}.

\noindent
{\it Outline of the proof:} Using the coarse-grained approximation $\BARH^{(0)}_M(\eta)$ the decomposition \eqref{general} 
can  be rewritten as $R(\eta)=e^{-\BARH^{(0)}_M(\eta)}$, and thus we obtain
\[
\BARH_M(\eta)=\BARH_M^{(0)}(\eta)-\frac{1}{\beta}\log\int A(\sigma)\nu(d\sigma|\eta)\COMMA
\]
where $A$, and $\nu_{\eta}$ are given abstractly in \eqref{abstractdecomp} and are defined both for short and long-range interactions 
in analogy to \eqref{nu}.
The construction of the series in \eqref{expansion} relies on the  cluster expansion of the type
\begin{equation}\label{decomp4}
A(\sigma)\equiv \prod_{i<j}(1+f_{ij}^l)\prod_{i:\,\rm{odd}}(1+f_{i-1 i+1}^s)=\sum_{G\in\mathcal{G}_M}\prod_{\{i,j\}\in E(G)}\tilde{f}_{ij}
\end{equation}
where
$$
\tilde{f}_{ij}=\begin{cases}
                 f_{ij}^l \,\,\text{or}\,\,\,f_{ij}^s,  & \text{if $i$ even and $j=i\pm 2$} \\  
                 f_{ij}^l & \text{otherwise,}
       \end{cases}
$$
and $\mathcal{G}_M$ is the set of all graphs on $M$ vertices, where $M$ is the total number of coarse cells.
Such an equality and the complete proof  is
carried out  in Section~\ref{polymer}. 
In turn, the terms on the right hand side of \eqref{decomp4} give rise to the  expansion \eqref{expansion} and the corresponding 
higher-order corrections. 

\subsection{ A posteriori error estimates} In \cite{KPRT} we introduced the use of cluster expansions as a
tool for constructing a posteriori error estimates for coarse-graining problems, based on the rather simple observation that
higher-order terms in \eqref{expansion} can be viewed as errors that depend only on the coarse variables $\eta$. 
Following the same approach an a posteriori estimate immediately follows from \eqref{expansion}.
\begin{corollary}
We have
\[
  \RELENTR(\BARIT{\mu}_{M,\beta}^{(0)}\SEP\mu_{N,\beta}\circ \COP^{-1})=
  \beta\EXPECT_{\BARIT{\mu}_{M,\beta}^{(0)}}[S(\eta)]+
  \log\left(\EXPECT_{\BARIT{\mu}_{M,\beta}^{(0)}}[\EXP{-\beta S(\eta)}]\right)
  +\BIGO(\delta^2)\COMMA
\]
where the residuum operator is  $S(\eta)=\BARH^{(1)}_M(\eta)$.

\end{corollary}
In  \cite{KPRT2} we already employed this type of estimates  for stochastic lattice systems
with long-range interactions, in order to construct adaptive coarse-graining schemes.  These tools  
operated as an ``on-the-fly'' coarsening/refinement method that recovers accurately phase-diagrams. 
The estimates allowed us to change adaptively the coarse-graining  level within the coarse-graining  hierarchy
once sufficiently large or small errors were detected, thus  speeding up the calculations of phase diagrams.
Earlier work that uses only an upper bound and not the asymptotically sharp cluster expansion-based estimate 
can be found in \cite{CKV1,CKV2}. 

\subsection{Microscopic reconstruction}
The reverse procedure of coarse-graining, i.e. reproducing ``atomistic" properties,  directly from coarse-grained  simulation methods is an issue that arises extensively  in the polymer science literature, \cite{tsop, mp}. The principal idea  is that computationally inexpensive coarse-graining algorithms   
will reproduce  large scale structures and subsequently microscopic information will  be added through 
{\em microscopic reconstruction},  for example
the calculation of  diffusion of penetrants through polymer melts, reconstructed from  CG simulation, \cite{mp}.

In this direction,  the CGMC methodology discussed in this section can provide a  framework to mathematically formulate microscopic  reconstruction  and study related numerical  and computational issues.
Indeed, the conditional measure $A(\sigma)\nu(d\sigma|\eta)$ in the  multi-scale decompositions \eqref{general} and \eqref{general-long}
can be also viewed as a microscopic reconstruction of the Gibbs state \eqref{microGibbs} once 
the coarse variables $\eta$ are specified.
The product structure in \eqref{decomp2}
and \eqref{decomp3} allows for easy generation of the fine scale details by first
reconstructing over a family of domains 
given only the coarse-grained data and
gradually moving to the next family of domains 
given now both the coarse-grained data and the previously reconstructed microscopic values. 

In view of of this abstract procedure 
based on  multiscale decompositions such as \eqref{general}, we readily see that the particular product  structure of the
explicit formulas \eqref{aa} and  \eqref{nu} for the case of the dimension $d=1$
yields  a hierarchy of reconstruction schemes.
 A  first order approximation can be based on the approximation  $A(\sigma)\sim 1$ (cf.  \eqref{asymptotics}, \eqref{aa}): 
 \smallskip
 
 \begin{itemize}
 \item[(a)] first,  $R(\eta)$ defined in \eqref{aa} provides  the  coarse-graining scheme, which will produce  coarse variable data $\eta(k)$ for all $k$;  
 
 \item[(b)] next, 
   we  reconstruct the microscopic configuration $\sigma^{\text{even}}$ consisting of  the $\sigma^k$'s in all boxes (coarse-cells) 
with $k$ {\it even} using  the measure 
$\nu_k(d\sigma|\eta):=\frac{e^{-\beta H^s_k}  Z_{k+1}(\sigma^k , 0) Z_{k-1}(0,\sigma^k)} {Z_{k-1,k,k+1}(0,0)}P_k(\sigma^k\SEP \eta(k))$, conditioned on the coarse configuration $\eta(k)$ from (a) above;

\item[(c)] finally,   we reconstruct the microscopic configuration in the remaining  boxes with
$k$ {\it odd} using  $\nu_k(d\sigma|\eta):=\frac{ e^{-\beta \big( H^s_k + W_{k-1,k} + W_{k,k+1}\big)}}{ Z_k(\sigma^{k-1},\sigma^{k+1})} 
P_k(d \sigma^k\SEP \eta(k))$, given the coarse variable $\eta(k)$ from step (a), and  the microscopic configurations $\sigma^{\text{even}}$ from step (b).

\end{itemize}
\smallskip

\noindent
We note that this procedure is {\em local} in the sense that the reconstruction can be carried out in only the  ``subdomain  of interest" of the 
entire microscopic lattice $\LATT$; this is  clearly  computationally advantageous because microscopic kMC solvers are used only in the specific   part of the computational domain, while inexpensive CGMC solvers are used in the entire coarse lattice $\LATTC$.

Further discussion on the numerical analysis issues related  to microscopic reconstruction for lattice systems with long-range interactions can be found in 
\cite{KPR, TT, KPS, KT}.

\section{Semi-analytical coarse-graining schemes and examples}\label{semi}
Next  we discuss  the numerical implementation  of the 
effective coarse-grained Hamiltonians derived in Theorem~\ref{maintheorem}.
We begin with a general implementation scheme and we subsequently
investigate further simplifications for particular examples in one 
space dimension.

\subsection{ Semi-analytical splitting schemes and inverse Monte Carlo methods} 
One of the main points of our method is encapsulated in (\ref{split}): the
computationally expensive long-range part for conventional Monte Carlo methods
can be computed by calculating the analytical formula given in
\eqref{l01} in the spirit of our previous work \cite{KPRT}.
Then we can turn our attention to the short-range interactions where Monte Carlo methods, 
at least for reasonably sized domains, are inexpensive. More specifically for the evaluation 
of the short-range contribution in \eqref{split} we introduce the normalized measure
\begin{equation}\label{Phat}
\hat{P}_k(d\sigma^k\SEP\eta(k))=\frac{1}{Z_k(\eta(k); 0, 0)}e^{-\beta H^s_{k}}
P_k(d\sigma^k\SEP\eta(k))\COMMA
\end{equation}
where the sum is computed with free boundary conditions on $\CUBE_k$ and $Z_k(\eta(k); 0, 0)$ 
is accordingly defined as in (\ref{z0boundary}). Thus (\ref{decomp}) can be rewritten as 
\begin{equation}\label{semi0}
\bar{H}_M^{s,(0)}=\sum_{k\in\bar{\Lambda}}\bar{U}_k^{s,(0)}(\eta(k))+\sum_{k:\,\mathrm{even}}
\bar{V}_{k-1,k,k+1}^{s,(0)}(\eta(k-1), \eta(k),  \eta(k+1))\COMMA  
\end{equation}
where, based on (\ref{decomp}) and (\ref{Phat}), we defined the  three-body coarse interaction potential
\begin{eqnarray}
\bar{V}_{k-1,k,k+1}^{s,(0)}&&(\eta(k-1),\eta(k),  \eta(k+1))= -\frac{1}{\beta}\log\int e^{-\beta (W_{k-1,k}(\sigma)+W_{k,k+1}(\sigma))}
\nonumber
\\
&&\times\hat{P}_{k-1}(d\sigma^{k-1}\SEP\eta(k-1))
\hat{P}_k(d\sigma^k\SEP\eta(k))\hat{P}_{k+1}(d\sigma^{k+1}\SEP\eta(k+1))
 \PERIOD \label{semi1}
\end{eqnarray}
The main difficulty in the calculation of (\ref{semi1}) is that for the three-body integral one needs to perform the 
integration for all possible combinations of the multi-canonical constraint. On the other hand all simulations involve 
only short-range interactions and need to be carried out only on  three coarse cells, rather than the entire lattice. 
Practically, the calculation of (\ref{semi1}) can be implemented   using the so-called {\em inverse  Monte Carlo method}, 
\cite{kremerplathe}.
We sample the measure $\hat{P}_{k}$ using  Metropolis spin flips and subsequently 
we create a histogram for all possible values of $\eta(k)=\sum_{x\in\CUBE_k}\sigma(x)$. Then
we compute the above integral by using the samples which correspond to the
prescribed values $\eta(k-1),\eta(k)$ and $\eta(k+1)$.

A complementary approach in order to further increase the computational efficiency of the schemes presented in Theorem~\ref{maintheorem} 
is to rearrange the splitting  based on the size of the error  in \eqref{relent}. Indeed, these estimates
suggest a natural way to decompose the overall interaction potential
into: (a)  a short-range piece $J_s$ including  possible singularities originally in $J$, e.g., the non-smooth part in the Lennard-Jones potential,
and (b) a locally integrable (or smooth) long-range decaying component, $J_l$. Thus, if $K(x, y)$ is the short-range
potential in \eqref{microHamiltonian} we can rewrite the overall potential as
\begin{equation}
\label{split2}
K(x,y)+J(x,y) = J_{s}(x,y) + J_{l}(x,y)\PERIOD
\end{equation}
In this way the accuracy can be enhanced by implementing the analytical coarse-graining \eqref{l02} for the smooth  
long-range piece  $J_{l}(x,y)$, and the semi-analytical scheme \eqref{decomp} for the ``effective'' short-range piece  $J_{s}(x,y)$. 

\begin{remark}
{\rm
Existing methods, e.g., \cite{doi}, employ an inverse Monte Carlo computation involving both short and long-range interactions, 
and due to computational limitations have to disregard multi-body terms such as the ones considered in the method proposed here.
The splitting approach developed here allows us to calculate analytically 
the approximate effective Hamiltonian for the costly long-range interactions, \VIZ{l02} in \VIZ{split} or \eqref{split2}, 
and in parallel carry out the inverse Monte Carlo  step for  \VIZ{semi0}.
The necessity to  include multi-body terms in the effective Hamiltonian  was first discussed in \cite{AKPR} together with their  
role in the proper coarse-graining of singular short-range interactions. We further quantify the regimes where such
multi-body  terms are necessary in the context of a specific example.
}
\end{remark}

\subsection{A typical example: improved schemes  and a posteriori estimation}
We examine  the derived coarse-graining schemes in the context of a specific, but rather typical example.
We consider the Hamiltonian
\begin{equation}\label{nnH}
H_N(\sigma)=H^{s}_N(\sigma)+H^{l}_N(\sigma):=K\sum_{\langle x,y\rangle}\sigma(x)\sigma(y)-\frac{1}{2}\sum_{(x,y)}J(x-y)\sigma(x)\sigma(y)
\end{equation}
where by $\langle x,y\rangle$ we denote summation over the nearest neighbors, i.e., $|x-y|=1$,  and by $(x,y)$ 
the long range summation as in \eqref{defJV1}. 
Although we follow the splitting strategy discussed in the previous paragraph we present a  simplified numerical algorithm 
by carrying out further analytical calculations. Not surprisingly, such calculations allow not only for easier sampling in 
the semi-analytical calculations of the inverse Monte Carlo, but give additional insight on the nature of multi-body, 
coarse-grained interactions.

For the short-range contributions,
given a coarse cell $\CUBE_k$ with $q$ lattice points, we denote by $x_1,\ldots, x_q$ the
lattice sites in $\CUBE_k$. With this notation, following \eqref{semi1} the short-range three-body interaction is given by
\begin{eqnarray}
\bar{V}_{k-1,k,k+1}^{s,(0)}&&(\eta(k-1),\eta(k),  \eta(k+1))= -\frac{1}{\beta}\log\int e^{-\beta K(\sigma^{k-1}(x_q)\sigma^k(x_1)+\sigma^k(x_q)\sigma^{k+1}(x_1))}
\nonumber
\\
&&\times\hat{P}_{k-1}(d\sigma^{k-1}\SEP\eta(k-1))
\hat{P}_k(d\sigma^k\SEP\eta(k))\hat{P}_{k+1}(d\sigma^{k+1}\SEP\eta(k+1))
 \PERIOD \label{semi2}
\end{eqnarray}

The main difficulty in computing the second term is the conditioning on  the coarse-grained values
$\eta(k-1), \eta(k), \eta(k+1)$ over three coarse cells.
At first glance this requires to run multi-constrained Monte Carlo 
dynamics for every given value of the $\eta$'s, i.e., for $q^3$ variables.
However, as we show in the sequel,
when dealing with a particular example, e.g., the nearest neighbor interactions, the computationally
expensive three-body term reduces to product of one-body terms. We first rewrite
\[
e^{-\beta K\sigma^{k-1}(x_q)\sigma^k(x_1)}=a-b\sigma^{k-1}(x_q)\sigma^k(x_1)\COMMA
\]
where we set
\[
\label{lambda}a=\cosh(\beta K)\, , \quad b=\sinh(\beta K)\, , \quad 
\lambda=\tanh(\beta K)\PERIOD
\]
Moreover, we introduce the one- and two-point correlation functions
\[
\Phi^x_k(\eta_k):=\int \sigma(x)\hat{P}_{k}(d\sigma^k\SEP\eta(k))
\,\,\,\,\,\,\,\mathrm{and}\,\,\,\,\,\,\,
\Phi^{x,y}_k(\eta_k):=\int \sigma(x)\sigma(y)\hat{P}_{k}(d\sigma^k\SEP\eta(k))\, .
\]
By symmetry we have that $\Phi^{x_1}_k=\Phi^{x_q}_k$ and 
similarly,  consider $\Phi^{x_1,x_q}_k$ for  $x=x_1$ and $y=x_q$. 
Furthermore, these functions depend on $k$ only via the coarse variable $\eta_k$,
hence we now define
\begin{equation}\label{phi}
\Phi^1(\eta_k):=\int \sigma(x_1)\hat{P}_{k}(d\sigma^k\SEP\eta(k))
\,\,\,\,\,\,\,\mathrm{and}\,\,\,\,\,\,\,
\Phi^{2}(\eta_k):=\int \sigma(x_1)\sigma(x_q)\hat{P}_{k}(d\sigma^k\SEP\eta(k))\PERIOD
\end{equation}
It is a straightforward computation to show that 
\begin{eqnarray}
\bar{V}_{k-1,k,k+1}^{s,(0)}(\eta(k-1), &&\eta(k), \eta(k+1))
=
-\frac{2}{\beta}\log a - \nonumber
\\
-\frac{1}{\beta}\log \Big(&&
1 - \lambda \Phi^1(\eta(k-1))\Phi^1(\eta(k))
-\lambda \Phi^1(\eta(k))\Phi^1(\eta(k+1))\nonumber
\\
&&+\lambda^2 \Phi^1(\eta(k-1))\Phi^2(\eta(k))\Phi^1(\eta(k+1))
\Big)\label{newformula}
\end{eqnarray}
Although these are three-body interactions, the additional analytical calculations reduce their computation 
to the nearest-neighbor Monte Carlo 
sub-grid sampling  of \eqref{phi}.
Moreover, from \eqref{smallness} we have
\[
f^s_{k-1,k+1}(\eta(k); \sigma^{k-1}, \sigma^{k+1})  
=\frac{\lambda^2\sigma^{k-1}(x_q)\sigma^{k+1}(x_1)[\Phi^2(\eta(k)-(\Phi^1(\eta(k))^2]}
{(1-\lambda\sigma^{k-1}(x_q)\Phi^1(\eta(k)))
(1-\lambda\sigma^{k+1}(x_1)\Phi^1(\eta(k)))}
\]
thus the following estimate holds for some $C>0$
\begin{equation}\label{conditional}
\sup_{\sigma^{k-1},\sigma^{k+1}}|f^s_{k-1,k+1}|\leq C\lambda^2|\Phi^2(\eta(k))-[\Phi^1(\eta(k))]^2|\equiv\Theta(\eta_k; \lambda)\COMMA
\end{equation}
where the right-hand side $\Theta$ is an {\em a posteriori} functional in the sense that it
can be computed from the coarse-grained data.
In fact, we can  estimate the a posteriori error indicator by an analytical formula.
A high temperature expansion yields
\begin{eqnarray}
\Phi^1 (\eta_k)   &=&   \EXPECT[\sigma(x)\SEP\eta]+\BIGO(\lambda)=\frac{\eta}{q}+ \BIGO(\lambda) \label{Phi1} 
\\\nonumber
\\
\Phi^2 (\eta_k) &=& \EXPECT[\sigma(x)\sigma(y)\SEP\eta]+ \BIGO(\lambda)= \frac{\eta^2 - q}{q(q-1)}+ \BIGO(\lambda) 
\PERIOD \label{Phi2} 
\end{eqnarray}
Then,
\begin{equation}\label{conditional2}
\Theta(\eta_k; \lambda)\sim \lambda^2|\Phi^2_k-(\Phi^1_k)^2|=
\lambda^2{q^2-\eta^2\over q^2(q-1)}+\BIGO(\lambda^3)\, .
\end{equation}
Thus 
the validity of Theorem~\ref{maintheorem} and the  derived coarse-grained approximations can be {\em conditionally} 
checked during simulation by
\begin{equation}\label{conditional3}
\sup_{\sigma^{k-1},\sigma^{k+1}}|f^s_{k-1,k+1}|\leq C {\lambda^2\over q-1}\Big(1-{\eta^2\over q^2}\Big)+\BIGO(\lambda^3)\PERIOD
\end{equation}

We note that (\ref{conditional3}) suggests a quantitative understanding of the dependence of   
the coarse-graining error for the nearest-neighbor Ising model.  The error increases,  (a) 
 when the parameter $\lambda^2$ increases, i.e., at lower temperatures/stronger  short-range interactions, 
(b)  when the level of coarse-graining $q$ decreases, and
(c) at regimes where the local coverage $\eta$ is not uniformly homogeneous, i.e.,
away from the regime $\eta \approx \pm q$. Such situation occurs, for example, around an  interface in the phase transition regime. 
This is the case even in one dimension if long-range interactions are present in the system.

\section{Proofs}\label{polymer}

In this section we first construct and  prove the convergence of the cluster expansion. We formulate the
proofs in the full generality assuming a $d$-dimensional lattice. Thus coordinates of lattice
points are understood as multi-indices in $\Z^d$.
We start by constructing the a priori coarse-grained measure induced by the short-range interaction.
We perform a block decimation procedure following the strategy in \cite{op} and
partition $\LATTC$ into $2^d$-many sublattices of spacing $2q$.
\begin{figure}[h]
  \centerline{\hbox{\psfig{figure=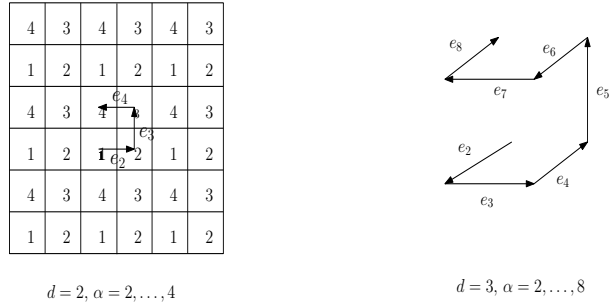,height=4cm,width=8cm}}}
  \caption{\label{cube} The sublattices $\LATTC^{\alpha}$ covering the coarse lattice $\LATTC$. The vectors $e_\alpha$ defining
   translations of the first sublattice $\LATTC^{1}$ are depicted for $d=2,3$. The cells on the two-dimensional lattice are numbered
   with values of $\alpha=1,\dots,4$ according to what sublattice $\LATTC^{\alpha}$ they belong.}
\end{figure}
Let $e_{\alpha}$, $\alpha=2,3,\ldots,2^d$ be vectors (of length $q$) along the edges of $\LATTC$ as demonstrated in 
Figure~\ref{cube} for $d=3$.
We write the coarse lattice as union of sub-lattices
\begin{equation}\label{partition}
\LATTC=\cup_{\alpha=1}^{2^d}\LATTC^{\alpha}\COMMA
\end{equation}
where $\LATTC^1=2\LATTC$, $\LATTC^2=\LATTC^1+e_2$ and $\LATTC^{\alpha+1}=\LATTC^{\alpha}+e_{\alpha+1}$, for $\alpha = 1,\dots,2^d-1$.
Given a coarse cell $\CUBE_k$ we define the set of neighboring cells by
\[
\partial\CUBE_k:=\cup_{\{l:\,\|l-k\|=1\}}\CUBE_l\COMMA
\]
where $\|l-k\|:=\max_{i=1,\ldots,d}|l_i-k_i|$. We also let $D_k:=\CUBE_k\cup\partial\CUBE_k$.

Given a sublattice $\LATTC^{\alpha}$ we denote by $\sigma^{\alpha}$ the microscopic configuration
in all the cells $\CUBE_k\in\LATTC^{\alpha}$ and by $\sigma^{>\alpha}$ the configuration in $\LATTC^{\beta}$ for all $\beta>\alpha$.
We also define a function $p:\LATTC\to\{1,\ldots,2^d\}$ 
such that for $k\in\LATTC$, we have $p(k)=\alpha$ if $\CUBE_k\in\LATTC^{\alpha}$.

We split the short-range part of \eqref{microHamiltonian}
\[
H_N^s(\sigma)=\sum_{\alpha}\sum_{k\in\LATTC^{\alpha}}H^s_k(\sigma^{\alpha})+
\sum_{\alpha}\sum_{k\in\LATTC^{\alpha}}W_k(\sigma^{\alpha};\,\sigma^{>\alpha})\COMMA
\]
where, for $k\in\LATTC^{\alpha}$, the  terms $H_k(\sigma^{\alpha})$ are the self energy on
the boxes $\CUBE_k$ given by
\[
H^s_k(\sigma^{\alpha})=\sum_{X\subset\CUBE_k}U_X(\sigma^{\alpha})\PERIOD
\]
Moreover, the energy due to the interaction of $\CUBE_k$ with the neighboring cells is given by
\[
W_k(\sigma^{\alpha};\,\sigma^{>\alpha})
=\sum_{X\subset D_k}U_X(\sigma^{\alpha}\vee\sigma^{>\alpha})\COMMA
\]
where $\sigma^{\alpha}\vee\sigma^{>\alpha}$ is the concatenation on 
$\LATTC^{\alpha}$ and $\LATTC^{>\alpha}$.
Now we construct the reference conditional measure $\nu_{\eta}$ under the constraint of a fixed averaged 
value $\eta=\{\eta(k)\}_{k\in\LATTC}$ on the coarse cells.

\medskip
\noindent
\underline{\it Step 1.}
The starting point is a product measure on $\CUBE_k$ for $k\in\LATTC^1$.
We let  $A_1(k)\equiv\CUBE_k$ and after appropriate normalization we obtain
\begin{eqnarray}\label{start}
e^{-H_N^s(\sigma)}\prod_{k\in\LATTC}P_k(d\sigma) & = &
\prod_{\alpha\geq 2}\prod_{k\in\LATTC^{\alpha}}\left(
e^{-H^s_k(\sigma^{\alpha})}e^{-W_k(\sigma^{\alpha};\sigma^{>\alpha})}
P_k(d\sigma^{\alpha}) \right)\times\nonumber\\
&&\prod_{k\in\LATTC^1}Z(A_1(k);\sigma^{>1};\eta(k))\,\,\,\nu^1_{>1}(d\sigma^1)
\end{eqnarray}
where
\begin{equation}\label{nu1}
\nu^1_{>1}(d\sigma^1):=\prod_{k\in\LATTC^1}\left[
\frac{1}{Z(A_1(k);\sigma^{>1};\eta)}
e^{-W_k(\sigma^1;\,\sigma^{>1})}
e^{-H_k(\sigma^1)}P_k(d\sigma^1)\right]
\end{equation}
is the new prior measure on $\LATTC^1$
with boundary conditions $\sigma^{>1}$ and the canonical constraint $\eta(k)$, $k\in\LATTC^1$.
The partition function
$$
Z(A_1(k);\sigma^{>1};\eta(k))=\int e^{-H^s_k(\sigma^1)}e^{-W_k(\sigma^1;\,\sigma^{>1})}P_k(d\sigma^1)
$$
depending on the boundary conditions $\sigma^{>1}$ on the set $\partial A_1(k)$
couples the configurations in $\CUBE_l$ with $l\in\partial A_1(k)$.
In particular, it couples the configurations $\sigma^2$ and gives rise to a new interaction
between them for which it will be shown that it is small due to the condition~\ref{condition}.

\medskip
\noindent
\underline{\it Step 2.}
Moving along the vector $e_2$ we seek the measure $\nu^2_{>2}$ on 
$\{+1,-1\}^{\cup_{k\in\LATTC^2}\CUBE_k}$.
Given the partition function $Z(A_1(k);\sigma^{>1};\eta(k))$ we denote by
$S^+_{k,e_2}Z$ the partition function 
on the same domain $A_1(k)$ as $Z$, but with new boundary conditions which are the same
as $Z$ in the $+e_2$ direction, free in the $-e_2$ and
unchanged in all the other directions.
Similarly, we denote by $S^-_{k,e_2}Z$ the partition function with free boundary conditions in the
direction $+e_2$ and by $S^0_{k,e_2}Z$ with free boundary conditions in both $\pm e_2$ directions.
With these definitions we have the identity
\begin{equation}\label{id}
Z(A_1(k);\sigma^{>1};\eta(k))=\frac{(S^+_{k,e_2}Z) (S^-_{k,e_2}Z)}{(S^0_{k,e_2}Z)}(1+\Phi^1_k)\COMMA
\end{equation}
where we have introduced the function $\Phi^1_k$ which contains the interaction
between the variables $\sigma^{>1}$, and  it is given by
\[
\Phi^1_k:=\frac{Z(A_1(k);\sigma^{>1};\eta(k)) (S^0_{k,e_2}Z)}{(S^+_{k,e_2}Z) (S^-_{k,e_2}Z)}-1\PERIOD
\]
In this way we split the partition function $Z$ into a part where the interaction
between the cells $\CUBE_{k-e_2}$ and $\CUBE_{k+e_2}$ is decoupled and an error part
which is to be small.
The terms in the second product contain all possible interactions in the set 
\begin{equation}\label{A2}
A_2(k)=\CUBE_{k-e_2}\cup \CUBE_k\cup\CUBE_{k+e_2}
\end{equation}
for $k\in\LATTC^2$ with the corresponding
partition function being given by
\[
Z(A_2(k);\sigma^{>2};\eta(k))=\int e^{-H^s_k(\sigma^2)}e^{-W_k(\sigma^2;\,\sigma^{>2})}(S^+_{k-e_2,e_2}Z) (S^-_{k+e_2,e_2}Z)P_k(d\sigma^2)
\]
all due to the condition~\ref{condition}.

The next step is to index
the new partition functions $(S^+_{k,e_2}Z)$ and $(S^-_{k,e_2}Z)$ (which are functions of $\sigma^2$ indexed
by $k\in\LATTC^1$) with respect to $k\in\LATTC^2$.
We have
\[
\prod_{k\in\LATTC^1}(S^+_{k,e_2}Z) (S^-_{k,e_2}Z)=\prod_{k\in\LATTC^2}(S^+_{k-e_2,e_2}Z) (S^-_{k+e_2,e_2}Z)\PERIOD
\]
Then if we neglect for a moment the error term $(1+\Phi^1_k)$, in order to define $\nu^2_{>2}$ we have to deal with the following terms
\[
\prod_{k\in\LATTC^1}(S^0_{k,e_2}Z)^{-1}
\prod_{k\in\LATTC^2}\left[e^{-H^s_k(\sigma^2)}e^{-W(\sigma^2;\,\sigma^{>2})}(S^+_{k-e_2,e_2}Z) (S^-_{k+e_2,e_2}Z)P_k(d\sigma^2)
\right]\PERIOD
\]
The terms in the second product contain all possible interactions in the set 
$A_2(k)$, given in \eqref{A2} for $k\in\LATTC^2$ with the corresponding
partition function being given by
\[
Z(A_2(k);\sigma^{>2};\eta(k))=\int e^{-H^s_k(\sigma^2)}e^{-W_k(\sigma^2;\,\sigma^{>2})}(S^+_{k-e_2,e_2}Z) 
S^-_{k+e_2,e_2}Z)P_k(d\sigma^2)\PERIOD
\]
By normalizing with this function we obtain the measure
\begin{eqnarray}\label{nu2}
\nu^2_{>2}(d\sigma^2) & = & \prod_{k\in\LATTC^2}
\left[\frac{1}{Z(A_2(k);\sigma^{>2};\eta(k))}\times\right.\nonumber\\
&& \left.
e^{-H^s_k(\sigma^2)}e^{-W_k(\sigma^2;\,\sigma^{>2})}
(S^+_{k-e_2,e_2}Z) (S^-_{k+e_2,e_2}Z)P_k(d\sigma^2)\right]\PERIOD
\end{eqnarray}
Note that the factor $(S^0_{k,e_2}Z)^{-1}$ depends on $\eta$ as well as on $\sigma^{>2}$ and hence we will need 
to further split it when we define a measure on the variables
on which it depends.
Summarizing the first two steps we have obtained that the left hand side of \eqref{start}
is equal to
\[
\left[
\prod_{k\in\LATTC^2}Z(A_2(k);\sigma^{>2};\eta(k))
\prod_{k\in\LATTC^1}(S^0_{k,e_2}Z)^{-1}
\prod_{k\in\LATTC^1}(1+\Phi^1_k)
\right]
\nu^2_{>2}(d\sigma^2)
\nu^1_{>1}(d\sigma^1)\PERIOD
\]

If we are interested in the case $d=1$, this would be the final expression.
However, for higher dimensions we need to repeat the above steps.
We give one more step in order to obtain more intuition on the relevant terms and
then we give the final expression in agreement with the result in \cite{op}.
The proof of the general formula is done with a recurrence argument on the number of steps
and for the details we refer to \cite{op}.

\medskip
\noindent
\underline{\it Step 3.}
To proceed in the next step along direction $e_3$ we split
$Z(A_2(k);\sigma^{>2};\eta(k))$ (which couples the configurations in $\CUBE_k$ with $p(k)=3$) in the same fashion as before. 
We have
\[
Z(A_2(k);\sigma^{>2};\eta(k))=\frac{(S^+_{k,e_3}Z)(S^-_{k,e_3}Z)}{(S^0_{k,e_3}Z)}(\Phi^3_k+1)
\]
where
\[
\Phi^3_k:=\frac{Z(A_2(k);\sigma^{>2};\eta(k))(S^0_{k,e_3}Z)}{(S^+_{k,e_3}Z)(S^-_{k,e_3}Z)}-1\PERIOD
\]
We further change the indices in such a way that they are expressed with respect to $k\in\LATTC^3$
and then we glue the partition functions on $\CUBE_k$, $A_2(k-e_3)$ and $A_2(k+e_3)$.
We define
\[
A_3(k):=\CUBE_k\cup A_2(k-e_3)\cup A_2(k+e_3)\COMMA
\]
and
\[
Z(A_3(k);\sigma^{>3};\eta(k)):=\int e^{-H^s_k}e^{-W_k(\sigma^3;\sigma^{>3};\eta)}(S^+_{k-e_3,e_3}Z)(S^-_{k+e_3,e_3}Z)
P_k(d\sigma^3)\PERIOD
\]
The corresponding measure is
\begin{eqnarray}\label{nu3}
\nu^3_{>3}(d\sigma^3) & = & \prod_{k\in\LATTC^3}\left[\frac{1}{Z(A_3(k);\sigma^{>3};\eta(k))}\times\right.
\nonumber\\
&& \left.
e^{-H^s_k(\sigma^3)}e^{-W_k(\sigma^3;\,\sigma{>3})}
(S^+_{k-e_3,e_3}Z) (S^-_{k+e_3,e_3}Z)P_k(d\sigma^3)\right]\COMMA
\end{eqnarray}
and the left hand side of \eqref{start} is now equal to
\begin{eqnarray*}
&& \prod_{k\in\LATTC^4}\left[ e^{-H^s_k(\sigma^4)}e^{-W_k(\sigma^4;\sigma^{>4})}
P_k(d\sigma^4)\right]
\prod_{k\in\LATTC^4}Z(A_3(k);\sigma^{>3};\eta(k))
\prod_{k\in\LATTC^2}(S^0_{k,e_3}Z)^{-1}\times \\
&& \prod_{k\in\LATTC^1}(S^0_{k,e_2}Z)^{-1}
\prod_{k\in\LATTC^2}(1+\Phi_k^3)\prod_{k\in\LATTC^1}(1+\Phi_k^1)
\nu^3_{>3}(d\sigma^3)
\nu^2_{>2}(d\sigma^2)
\nu^1_{>1}(d\sigma^1)\PERIOD
\end{eqnarray*}
As in the previous steps we need to perform the usual actions on the partition function
$Z(A_3(k);\sigma^{>3};\eta(k))$ which will give rise to a new element $A_4(k)$ with $k\in\LATTC^4$
and new error terms $\Phi^4_k$ with $k\notin\LATTC^4$.
Furthermore, a similar splitting has also to occur for the factor $(S^0_{k,e_2}Z)^{-1}$
which also depends on $\sigma^4$, since the zero boundary condition involves only the direction $e_2$.
Related calculations will involve all the terms of similar origin as long as we move
to new sublattices $\LATTC^{\alpha}$, with $\alpha>4$, depending on the dimension.

\begin{figure}
  \centerline{\hbox{\psfig{figure=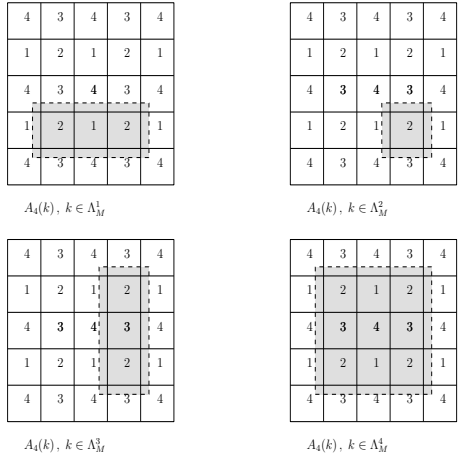,height=6cm,width=6cm}}}
  \caption{\label{A} The index sets $A_4(k)$ for $k\in\LATTC^{\alpha}$, $\alpha=1,2,3,4$, $d=2$, depicted
   as shaded cells. The cells in each lattice are numbered by $\alpha$ denoting the sublattice $\LATTC^{\alpha}$
   to which the cell belongs. }
\end{figure}

\smallskip 
\noindent{\sc Example: 2D lattice} The leading term in the approximation of
the coarse-grained Hamiltonian  $\BARH_M^{s}$ consists of terms that refer
to four different types of multi-cell interactions
\begin{eqnarray*}
\BARH_M^{s,(0)}&=&\sum_{k\in\LATTC^1}\log Z(A_4(k))-\sum_{k\in\LATTC^2}\log Z(A_4(k))\\
               &+&\sum_{k\in\LATTC^3}\log Z(A_4(k))-\sum_{k\in\LATTC^4}\log Z(A_4(k))\COMMA
\end{eqnarray*}
where $A_4(k)$ is a collection of coarse cells centered in $k\in\LATTC^{\alpha}$ 
and it is different depending on the sublattice to which the reference cell $k$ belongs. 
For $\alpha=1,2,3,4$ we have
$$
A_4(k) = \left\{\begin{array}{ll}
             \cup_{i,j\in\{-1,0,+1\}}\CUBE_{k+i e_2+j e_3}\COMMA & \;k\in\LATTC^4\COMMA \\
             \cup_{j\in\{-1,0,+1\}}\CUBE_{k+j e_3}\COMMA & \;k\in\LATTC^3\COMMA \\
             \CUBE_k\COMMA                               & \;k\in\LATTC^2\COMMA \\
             \cup_{i\in\{-1,0,+1\}}\CUBE_{k+i e_2}\COMMA & \; k\in\LATTC^1\PERIOD
                \end{array}\right.
$$
Figure~\ref{A} depicts the index sets $A_4(k)$ for the reference cell $k$ belonging
to $\LATTC^{\alpha}$ for $\alpha=1,\dots,4$. 

\medskip
\noindent
\underline{\it General formulation.}
At this point we proceed with the general formulation for given
$\alpha$ of the relevant quantities which
are the reference measure $\nu^{\alpha}_{>\alpha}(d\sigma^{\alpha})$,
the error term $\Phi^{\alpha}_k$, with $k\in\LATTC$
and the sets $A_{\alpha}(k)$ and $B_{\alpha}(k)$, with the latter being the 
relevant boundary of $A_{\alpha}$. The index $\alpha$ indicates the sublattice we are considering.
\begin{definition}\label{defAB}
The sets $A_{\alpha}(k)$ and $B_{\alpha}(k)$ for $k\in\LATTC^{\alpha}$ are
$$
A_{\alpha}(k) = \cup_{l: \|l-k\|=1,\, p(l)\leq\alpha}\CUBE_l\COMMA\;\;\;
B_{\alpha}(k) = \cup_{l: \|l-k\|=1,\, p(l)>\alpha}\CUBE_l\PERIOD
$$
\end{definition}

\begin{definition}\label{def_nugen}
Given $\alpha=1,\ldots,2^d$ we define the normalized Bernoulli measure on $\LATTC^{\alpha}$
\begin{equation}\label{nugen}
\nu^{\alpha}_{>\alpha}(d\sigma^{\alpha})=\prod_{k\in\LATTC^{\alpha}}\nu^{\alpha}_{B_{\alpha}(k)}(d\sigma^{\alpha})\COMMA
\end{equation}
where
\begin{equation}\label{nugenk}
\nu^{\alpha}_{B_{\alpha}(k)}(d\sigma^{\alpha})=
\frac{e^{-H^s_k(\sigma^{\alpha})}e^{-W_k(\sigma^{\alpha};\,\sigma^{>\alpha})}}{Z(A_{\alpha}(k);\sigma^{>\alpha};\eta(k))}
Z(A_{\alpha}(k)/\{k\};\sigma^{>\alpha};\eta(k))\prod_{l\in B_{\alpha}(k)}P_l(d\sigma^{\alpha})\PERIOD
\end{equation}
\end{definition}

As we have seen in {Step 3} we have two kinds of error terms $\Phi^{\alpha}_k$, in particular, those
with $k\in\LATTC^{\alpha}$ and others with $k\notin\LATTC^{\alpha}$.
In order to describe the latter we need to introduce additional notation.

For $\alpha=1,\ldots,2^d$ we denote by $\Gamma_{\alpha}$ the family of parallel hyperplanes
of dimension $d-1$ orthogonal to $e_{\alpha+1}$ passing through all the points $k\in\LATTC^{\alpha}$.
Note that for any $\alpha$,  we have that $\LATTC=\Gamma_{\alpha}\cup(\Gamma_{\alpha}+e_{\alpha+1})$.
In the next definition we introduce a new parameter $\epsilon_{\alpha}(k)\in\{\pm 1\}$
depending on whether we should perform gluing or unfolding as discussed before.
This is determined as follows: for fixed $\alpha\in\LATTC$ let $d(\alpha,\beta)$ be
the distance between the sublattices $\LATTC^{\alpha}$ and $\LATTC^{\beta}$
in the metric $\|\alpha-\beta\|_\infty=\sum_{i=1}^d|\alpha_i-\beta_i|$.
Moreover, we can find orthogonal vectors $\{v_j\}_{j=1,\ldots,d(\alpha,\beta)}$
and a family of signs $\{\epsilon_j\}_{j=1,\ldots,d(\alpha,\beta)}$ such that
\[
\LATTC^{\alpha}=\LATTC^{\beta}+\gamma(\alpha,\beta)\,\,\,\,\,\,
\mbox{with}\,\,\,\,\, \gamma(\alpha,\beta)=\sum_{j=1}^{d(\alpha,\beta)}\epsilon_j \,v_j\PERIOD
\]
Note also that $|\gamma(\alpha,\beta)|=d(\alpha,\beta)$.
Then the exponents $\epsilon_{\alpha}(k)$ with $p(k)=\beta$ are given by
\[
\epsilon_{\alpha}(k):=(-1)^{|\gamma(\alpha,\beta)|}\PERIOD
\]
Furthermore, we denote by $Y(k,\gamma(\alpha,\beta))$ the affine hyperplane of codimension 
$|\gamma(\alpha,\beta)|$ orthogonal to the connecting vectors 
$\{v_j\}_{j=1,\ldots,|\gamma(\alpha,\beta)|}$
and passing through the point $k$
\[
Y(k,\gamma(\alpha,\beta))=\cap_{j=1}^{|\gamma(\alpha,\beta)|}Y(k,v_j)\COMMA
\]
where $Y(k,v)$ is the hyperplane of dimension $d-1$ passing through $k$
and being perpendicular to the vector $v$.
From the set of coarse-lattice points belonging to
$Y(k,v)$ we  define the corresponding set by
\[
 \mathcal Y(k,\gamma(\alpha,\beta)):=\cup_{l\in Y(k,\gamma(\alpha,\beta))}\CUBE_l\PERIOD
\]
Then, letting $k\in\LATTC^{\alpha}$, for $l$ such that $\CUBE_l\subset\partial\CUBE_k$
and with $l\in\LATTC^{\beta}$, for some $\beta$, we define
\begin{eqnarray}
A_{\alpha}(l) &=& \begin{cases}
                 \emptyset  & \text{if $p(l)>2^{d(\alpha)}$} \\  
                 A_{\alpha}(k) \cap \mathcal Y(l,\gamma(\alpha,\beta)) & \mbox{otherwise,}
       \end{cases}
\\
B_{\alpha}(l) &=& B_{\alpha}(k)\cap \mathcal Y(l,\gamma(\alpha,\beta))\PERIOD
\end{eqnarray}
With the above definitions we can determine the error terms in the general expansion.
\begin{definition}\label{errors}
For any $k\in\LATTC^{\alpha}$ and for $k\in\Gamma_{\alpha}$ the error terms are given by
$$
\Phi^{\alpha}_k=-1+
\frac{Z(A_{\alpha}(k);\sigma^{>\alpha};\eta(k))Z(A_{\alpha+1}(k);\sigma^{>\alpha+1};\eta(k))}
{(S^+_{k,e_{\alpha+1}}Z)(S^-_{k,e_{\alpha+1}}Z)}\PERIOD
$$
Moreover, if $k\in\Gamma_{\alpha}+e_{\alpha+1}$ and $k\notin\LATTC^{\alpha+1}$ we have:
$$
\Phi^{\alpha}_k=-1+
\left[
\frac{Z(A_{\alpha}(k);\sigma^{>\alpha};\eta(k))Z(A_{\alpha+1}(k);\sigma^{>\alpha+1};\eta(k))}
{(S^+_{k-e_{\alpha+1},e_{\alpha+1}}Z)(S^-_{k+e_{\alpha+1},e_{\alpha+1}}Z)}
\right]^{-\epsilon_{\alpha}(k)}\PERIOD
$$
Furthermore, if $k\in\LATTC^{\alpha+1}$ we replace $Z(A_{\alpha+1}(k);\sigma^{>\alpha+1})$
by $Z(A_{\alpha+1}(k)/\{k\};\sigma^{>\alpha+1})$.
\end{definition}

From Proposition~2.5.1 in \cite{op} we have that 
the general $d$-dimensional formulation of the a priori measure induced by 
the short-range interactions is
\[
e^{-H_N^s(\sigma)}\prod_{k\in\LATTC}P_k(d\sigma)=
R^s(\eta)A(\sigma)\nu(d\sigma|\eta)\COMMA
\]
where we have the following factors
\begin{description}
\item[{\rm(i)}] a product of partition functions (depending only on the coarse-grained variable $\eta$) over
finite sets of coarse cells with supports $A_{2^d}(k)$, with $k\in\LATTC^{\alpha}$ and $\alpha=1,\ldots,2^d$
\begin{equation}\label{R}
R^s(\eta):=
\prod_{\alpha=1}^{2^d}
\prod_{k\in\LATTC^{\alpha}}
\left[
Z(A_{2^d}(k);\,\eta(k))^{\epsilon_{2^d}(k)}
\right]\COMMA
\end{equation}
\item[{\rm(ii)}] error terms in the form of a gas of polymers (with the only interaction to be a hard-core exclusion)
\[
A(\sigma):=\prod_{\alpha=1}^{2^d}\prod_{j\leq 2^{d(\alpha)}}\prod_{k\in\LATTC^j}(1+\Phi^{\alpha}_k)\COMMA
\]
\item[{\rm(iii)}] a reference measure induced by only the short-range interactions once we neglect the
reference system and the error terms
\[
\nu(d\sigma|\eta):=\nu^{2^d}\ldots\nu^2_{>2}\nu^1_{>1}\PERIOD
\]
\end{description}

With this expansion for the short range interactions, going back to the general strategy presented in Section~\ref{motivation}, 
if we also consider the long-range contribution from \eqref{l01}, we obtain
\begin{eqnarray*}
e^{-\beta\BARH_M(\eta)} & = & \int e^{-\beta H_N^l}e^{-\beta H_N^s}\prod_kP_k(d\sigma)\\
& = & e^{-\beta\BARH^{l,(0)}(\eta)}R(\eta)\int e^{-\beta (H_N^l-\BARH^{l,(0)})}A(\sigma)\nu(d\sigma|\eta)\COMMA
\end{eqnarray*}
which implies that
\begin{equation}\label{perturbation}
\BARH_M(\eta)=\BARH^{l,(0)}(\eta)
-\log R(\eta)
-\frac{1}{\beta}\log\EXPECT_{\nu}[ e^{-\beta(H_N^l-\BARH^{l,(0)})}A(\sigma)|\eta]\PERIOD
\end{equation}

\subsection{Cluster expansion and effective interactions}\label{sec-cluster}
The goal of this section is to expand the term 
$\EXPECT_{\nu}[\EXP{-\beta(H_N^l(\sigma)-\BARH^{l,(0)}(\eta))}A(\sigma)|\eta]$ in \eqref{perturbation}
into a convergent series using a cluster expansion.
By the construction given previously
the terms in $A(\sigma)$ are already in the form of a polymer gas with hard-core interactions only.
For the long-range part we first write the difference $H_N^l(\sigma)-\BARH^{l, (0)}(\eta)$ as
\begin{eqnarray}
&& H_N^l(\sigma)-\BARH^{l,(0)}(\eta)= \sum_{k\leq l} \Delta_{kl}J(\sigma)\COMMA
   \;\;\;\;\mbox{where} \nonumber \\
&& \Delta_{kl}J(\sigma):= -\HALF\sum_{\genfrac{}{}{0cm}{2}{x\in \CUBE_k}{y\in \CUBE_l,y\neq x}}
   (J(x-y)-\bar{J}(k,l))\sigma(x)\sigma(y)(2-\delta_{kl})\PERIOD \label{star}
\end{eqnarray}
We also define $f_{kl}^{}(\sigma):=\EXP{-\beta\Delta_{kl}J(\sigma)}-1$ and we obtain
\begin{equation}\label{first}
\EXPECT_{\nu}[\EXP{-\beta(H_N^l(\sigma)-\BARH^{l,(0)}(\eta)}A(\sigma)|\eta] =
\int\prod_{k\leq l}(1+f_{kl})
\prod_{\alpha=1}^{2^d}\prod_{j\leq 2^{d(\alpha)}}\prod_{k\in\LATTC^j}(1+\Phi^{\alpha}_k)
\nu(d\sigma|\eta) \PERIOD
\end{equation}

We define the polymer model which contains combined interactions originating from both
the short and long-range potential.
By expanding the products in \VIZ{first} we obtain terms of the type
\[
\prod_{j=1}^p\Phi^{\alpha_j}_{k_j}\prod_{i=1}^q f_{l_i,m_i}\,\,\,\,\,\,
\mbox{where}\,\,\, k_j,l_i,m_i\in\LATTC\,\,\,\mbox{and}\,\,\,\alpha_j\in\{1,\ldots,2^d\}
\]
for some $p$ and $q$.
The factors $\Phi^{\alpha_j}_{k_j}$ are functions of the variables which are on the boundary
of the corresponding sets $A_{\alpha_j}(k_j)$. This boundary is described by the set
\begin{equation}\label{C0}
C_0^{\alpha}(k)=\begin{cases}
                 B_{\alpha}(k)  & \text{if $k\in\Gamma_{\alpha}$,} \\  
                 B_{\alpha+1}(k) & \text{if $k\in\Gamma_{\alpha}+e_{\alpha+1}$.}
       \end{cases}
\end{equation}
Furthermore, since the measure $\nu(d\sigma|\eta)$ is not a product measure but
instead a composition of measures each one parametrized by variables which
are integrated by the next measure, 
we need to create a ``safety'' corridor around the sets $C_0^{\alpha}$
depending on the level of $\alpha$.
This is given in the next definition. For a given integer $\beta$,  with $1<\beta<2^d-\alpha$ we define
\begin{equation}\label{Cbeta}
C_{\beta}^{\alpha}(k)=
\cup_{\epsilon_1,\ldots,\epsilon_{\beta}\in\{\pm 1\}^{\beta}}
\cup
_{l:\CUBE_l\subset\partial (\CUBE_{k+\epsilon_1 e_{\alpha+1}+\ldots+\epsilon_{\beta}e_{\alpha+\beta}}),\, p(l)>\alpha+\beta}\CUBE_l\COMMA
\end{equation}
Then for given $\alpha\in\{1,\ldots,2^d\}$ we call a ``bond'' of type $C^{\alpha}$ the set
\begin{equation}\label{bond}
C^{\alpha}(k)=\cup_{\beta=0}^{2^d-\alpha}C_{\beta}^{\alpha}(k)\PERIOD
\end{equation}
With this definition, any factor $\Phi^{\alpha_j}_{k_j}$ has a region of dependence which is given
by the bond $C^{\alpha}(k)$.
Similarly, for the factors $f_{l_i,m_i}$ originating from the long-range interactions the initial
domain of dependence is $\CUBE_{l_i}\cup\CUBE_{m_i}$.
However, due to the non-product structure of the measure we need to introduce a safety corridor in the same way.
Given $k\in\LATTC$ for $\beta$ an integer with $1<\beta<2^d-p(k)$ we define
\begin{equation}\label{Ck}
C_{\beta}(k)=
\cup_{\epsilon_1,\ldots,\epsilon_{\beta}\in\{\pm 1\}^{\beta}}
\cup
_{l:\CUBE_l\subset\partial (\CUBE_{k+\epsilon_1 e_{\alpha+1}+\ldots+\epsilon_{\beta}e_{\alpha+\beta}}),\, p(l)>p(k)+\beta}\CUBE_l  \PERIOD
\end{equation}
Then for a given $f_{kl}$ we define
\begin{equation}\label{Ckl}
C(k,l)=\cup_{\beta=1}^{2^d-p(k)}C_{\beta}(k)\cup_{\beta=1}^{2^d-p(l)}C_{\beta}(l)\PERIOD
\end{equation}
With a slight abuse of notation we define for $R_0=\{k_1,\ldots,k_{|R_0|}\}$
\begin{equation}\label{CR}
C(R_0)=\cup_{i=1}^{|R_0|}\cup_{\beta=1}^{2^d-p(k_i)}C_{\beta}(k_i)\PERIOD
\end{equation}

A bond $l$ will be either a $C^{\alpha}_k$ bond for some $\alpha,k$, called of type $1$, 
or any subset of $\LATTC$, we call it a bond of type $2$.
We say that two bonds $l_1$ and $l_2$ are connected if $l_1\cap l_2\neq\emptyset$.
We call a polymer $R$ a set of bonds $l_1,\ldots,l_p,l_{p+1}$ where $l_1,\ldots,l_p$ are bonds of type $1$ and
$l_{p+1}$ is a bond of type $2$, i.e., a possibly empty subset $R_0\subset\LATTC$.
A polymer is called {\it connected} if
for all $i,j$, with $1\leq i<j\leq p+1$, there exists a chain of connected bonds in $R$ joining $l_i$
to $l_j$.
For such a polymer $R$ we define its cardinality 
to be two integers, the first counting the number of bonds of the type $1$ and the second
being the number of coarse cells included in the bond of type $2$, i.e.,
$\card(R):=(p,|R_0|)$.
The support $\supp(R)$ of $R$ is $\supp(R)=\cup_{i=1}^{p+1}l_i$
where $l_{p+1}\equiv C(R_0)$ (see \VIZ{CR}).
Let $\mathcal{R}$ be the set of all such polymers. Two polymers $R_1, R_2$ are said to be {\it compatible}
if $\tilde{R}_1\cap\tilde{R}_2=\emptyset$ and we write $R_1\sim R_2$. 

Given a polymer $R=l_1,\ldots,l_p,l_{p+1}$ we define the {\em activity} of $R$ to be the
function $w:\,\mathcal{R}\to\mathbb{C}$ given by
\begin{equation}\label{activity}
w(R)=\int \nu(d\sigma|\eta)
\left(\prod_{j=1}^p\Phi^{\alpha_j}_{k_j}
\sum_{g\in\mathcal{G}_{R_0}}\prod_{\{k,l\}\in E(g)} f_{k,l}\right)\COMMA
\end{equation}
where $\mathcal{G}_l$ is the collection of connected graphs on the vertices of $l$ (recall 
$l\subset\LATTC$) and $E(g)$ is the set of edges of the graph $g$.

We define a new graph $\mathbb{G}$ on $\mathcal{R}$ which has the edge $R_i$-$R_j$ if
the polymers 
$R_i$ and $R_j$ are not compatible.
We call $G\subset\mathcal{R}$ {\it completely disconnected} if the subgraph induced by $\mathbb{G}$
on $G$ has no edges. Let
\[
\mathcal{D}_{\mathcal{R}}=\cup_{n=0}^{|\mathcal{R}|}\{
(R_1,\ldots,R_n)\subset\mathcal{R}:\,\forall i\neq j,\, R_i \sim R_j\COMMA
\}
\]
then the partition function $Z$ can be written as
\[
Z=\sum_{G\in\mathcal{D}_{\mathcal{R}}}\prod_{R\in G}w(R)\COMMA
\]
which is the abstract form of a polymer model. Thus we can apply the general theorem of
the cluster expansion once we check the convergence condition.
The condition is stated as a theorem in \cite{BZ2}.
\begin{theorem}[\cite{BZ2}]
Let $a:\mathcal{R}\to\mathbb{R}_+$.
Consider the subset of $\mathbb{C}^{\mathcal{R}}$
\begin{eqnarray*}
\mathcal{P}_{\mathcal{R}}^a & := & \{
w(R),\, R\in\mathcal{R}:\,\forall \, R\in\mathcal{R}: \,\, |w(R)|e^{a(R)}<1\,\rm{and}\\
&& \sum_{R'\nsim R}(-\log(1-|w(R')|e^{a(R')}))\leq a(R)
\}\PERIOD
\end{eqnarray*}
Then on $\mathcal{P}_{\mathcal{R}}^a$, $\log Z$ is well defined and analytic and
\[
\log Z=\sum_{I\in\mathcal{I}(\mathcal{R})}c_{I}\prod_{R\in\rm{supp}(I)}w(R)^{I_R}\COMMA
\]
where $I=(I_R)_{R\in\mathcal{R}}$, $\mathcal{I}(\mathcal{R})$ is the collection of all multi-indexes $I$,
i.e., integer valued functions on $\mathcal{R}$, and
\[
c_I=\frac{1}{I_{R_1}!\ldots I_{R_{|\mathcal{R}|}}!}
\frac{\partial^{I_{R_1}+\ldots+I_{R_{|\mathcal{R}|}}}\log Z}{\partial^{I_{R_1}}w(R_1)\ldots\partial^{I_{R_{|\mathcal{R}|}}}
w(R_{|\mathcal{R}|})}|_{\{w(R_i)=0\}_i}\PERIOD
\]
\end{theorem}
For the proof we refer to \cite{BZ2}.
Thus we need to check the condition of convergence.
The following estimate for the long-range potential was proved in \cite{KPRT}.
\begin{lemma}\label{f}
Assume that $J$ satisfies \VIZ{defJV1}. Then
there exists a constant $C_1\sim\frac{q^{d+1}}{L}\|\GRADV\|_\infty$ such that
\begin{equation}\label{C}
\sup_{k\in\LATTC}\sum_{l:\,l\neq k}|\Delta_{kl}J(\sigma)|\leq C_1\COMMA
\end{equation}
for every $\sigma$.
\end{lemma}

For the short-range interaction we follow the analysis of \cite{op} and we consider the following condition
\begin{condition}\label{condition}
Let $e$ be a vector in one of the directions of the lattice $\LATTC$ and 
$Z_U(\Lambda;\sigma_-,\sigma_+,\tau;\eta_V)$ be the partition function 
for the interaction $U$ in the space domain $\Lambda$.
We consider boundary conditions $\sigma_{\pm}$ in the directions $\pm e$ and $\tau$ in all other
directions.
Moreover, we impose multi-canonical constraints $\eta(k)$ for $k\in V\subset\LATTC$ with
$\Lambda=\cup_{k\in V}\CUBE_k$.
For a given $q>r_0$, with $|\CUBE_k|=q^d$, the following inequality holds
\[
\sup_{\sigma_{\pm},\tau}
\sup_{\Lambda}
\sup_{\eta_V}
\Big|\frac{
Z_U(\Lambda;\sigma_-,\sigma_+,\tau;\eta_V)Z_U(\Lambda;0,0,\tau;\eta_V)}
{Z_U(\Lambda;0,\sigma_+,\tau;\eta_V)Z_U(\Lambda;\sigma_-,0,\tau;\eta_V)}-1\Big|
\leq C_2\COMMA
\]
where given the numbers $r=2^{2d}[3(2^{d+1}+1)]^d$, $E=2^{d+1}+1$, and $c>0$
the upper bound $C_2$ satisfies 
\[
r C_2 e^{c E}<1\PERIOD
\]
\end{condition}

Notice that we work with the same condition as {Condition $C_L$} defined in \cite{op}, 
where in our notation $L$ is $q$, yet similar analysis applies in order to prove convergence 
of the cluster expansion under the milder condition {\it Condition $C_L'$} again as in \cite{op}. 
We skip  the analysis of such issues since it goes beyond the goal of the present work.
Furthermore, these conditions are related to the ones presented in \cite{ds3} in order to ensure that
a given system belongs to the class of completely analytical interactions.
For further details we refer the reader to \cite{olivieri} and \cite{bco}
and to the references therein.

Now we are ready to prove the convergence condition.
\begin{lemma}\label{convergence}
The set $\mathcal{P}_{\mathcal{R}}^a$ is nonempty.
\end{lemma}

\PROOF
We take $a(R)=c |\supp(R)|$, where $c$ is a constant to be chosen later.
Note that $-\log(1-x)\leq 2x$, so it suffices to show that
\[
\sum_{R'\nsim R}2|w(R')| e^{a(R')}\leq a(R)\PERIOD
\]
Suppose that the generic polymer $R'$ is given by $R'=l_1,\ldots,l_p,l_{p+1}$,
for some $p\geq 0$,
where $l_j\equiv C^{\alpha_j'}(k_j')$ for $j=1,\ldots,p$
and $l_{p+1}=R_0'$, with $|R_0'|=n$ for some $n\geq 0$. 
For $|w(R')|$ we have
\[
|w(R')|\leq\int\nu(d\sigma|\eta)\prod_{j=1}^p|\Phi^{\alpha_j}_{k_j}|\cdot
|\sum_{g\in\mathcal{G}_{R_0'}}\prod_{\{k,l\}\in E(g)}f_{kl}|\PERIOD
\]
By the graph-tree inequality we have that for all $\sigma$, $\eta$ and with $|R_0'|=n$
\[
|\sum_{g\in\mathcal{G}_{R_0'}}\prod_{\{k,l\}\in E(g)}f_{kl}|
\leq \beta 2^n e^{n C_1}\sum_{\tau^0\in\mathcal{T}^0_n}
\prod_{\{k,l\}\in\tau^0}|\Delta_{kl}J(\sigma)|\COMMA
\]
where from Lemma~\ref{f} we have that $C_1\sim\frac{q^{d+1}}{L}\|\GRADV\|_{\infty}$.
We also let $\sup_{\sigma}|\Delta_{kl}J(\sigma)|\leq\Delta_{kl}$ with
$\Delta_{kl}\equiv q^{2d}\frac{1}{L^d}\frac{q}{L}\|\GRADV\|_{\infty}1_{(k,l):\,|l-k|\leq\frac{L}{q}}$.
Moreover, from Condition~\ref{condition} we have
\[
\int\nu(d\sigma|\eta)\prod_{j=1}^p|\Phi^{\alpha_j}_{k_j}|
\leq (C_2)^p\PERIOD
\]
Then for the activity $w(R')$ we obtain
\[
|w(R')|\leq \beta 2^n e^{n C_1}\sum_{\tau^0\in\mathcal{T}^0_n}
\prod_{\{k,l\}\in\tau^0}\Delta_{kl}\cdot (C_2)^p\PERIOD
\]
Thus to satisfy the sufficient condition for the convergence of the cluster expansion we
first bound the sum $\sum_{R'\nsim R}$ by
\[
\sup_{k\in R'}|\supp(R')|\sum_{p\geq 0}\sum_{n\geq 0}
\sum_{\genfrac{}{}{0cm}{2}{R':\,\supp(R')\supset\{k\},}{\card(R')=(p,n)}}\PERIOD
\]
The fixed coarse cell $\CUBE_k$ may belong to one of the $l_j$'s for $j=1,\ldots,p$
or to $R_0'$.
In the first case we estimate the sum over $R'$ by
\begin{eqnarray*}
&&\sum_{\genfrac{}{}{0cm}{2}{l_1:\,l_1\supset\{k\}}{l_2,\ldots,l_p}}\frac{1}{(p-1)!}
\sum_{\genfrac{}{}{0cm}{2}{k_1\in\,(\cup_j l_j)\cap R_0'}{k_2,\ldots,k_n}}\frac{1}{(n-1)!}\COMMA \\
&&\mbox{and in the second by}\\
&&\sum_{\genfrac{}{}{0cm}{2}{k_1=k}{k_2,\ldots,k_n}}\frac{1}{(n-1)!}
\sum_{\genfrac{}{}{0cm}{2}{l_1:\,l_1\cap R_0'\neq\emptyset}{l_2,\ldots,l_p}}\frac{1}{(p-1)!}\PERIOD
\end{eqnarray*}
Next, for every tree $\tau^0$ we have that
\[
\sup_{k\in R'}\sum_{\genfrac{}{}{0cm}{2}{k_1=k}{k_2,\ldots,k_n}}
\prod_{i,j\in\tau^0}\Delta_{k_i,k_j}\leq C_1^{n-1}\PERIOD
\]
We use the Cayley formula $\sum_{\tau^0\in\mathcal{T}^0_n}1=n^{n-2}$
and the fact that the cardinality of the sum 
$\sum_{\genfrac{}{}{0cm}{2}{l_1:\,l_1\supset\{k\}}{l_2,\ldots,l_p}}$
can be bounded by $r^p$, where $r=2^{2d}[3(2^{d+1}+1)]^d$ is an upper bound for the maximum number 
of $C^{\alpha}$ bonds that can pass through a point, as showed in \cite{op}.
Taking into account all the above we obtain
\begin{eqnarray*}
\sum_{R'\nsim R}2|w(R')| e^{a(R')} & \leq & |\supp(R')| \sum_{n\geq 1} A_1(n)
\sum_{p\geq 1} A_2(p)\COMMA
\end{eqnarray*}
where
\[
A_1(n)=\frac{1}{(n-1)!}n^{n-2}(\beta C_1)^{n-1}2^n e^{n\beta C_1}e^{c n}\,\,\,\,\,\,\,
\mathrm{and}\,\,\,\,\,
A_2(p)=\frac{1}{(p-1)!}r^p C_2^p e^{c E p}\COMMA
\]
and $E$ is the upper bound for the cardinality of any bond $C^{\alpha}(k)$, i.e.,
\[
\sup_{\alpha=1,\ldots,2^d}\sup_{k\in\LATTC^{\alpha}} |C^{\alpha}(k)|\leq E\equiv 2^{d+1}+1\PERIOD
\]
For $C_1$ and $C_2$ sufficiently small the two series converge to a finite number and we
choose $c$ to be this number.
\qed

\bibliographystyle{plain}
\bibliography{m1bib}

\end{document}

%% file: myppmacro3.tex
\def\RANGE{L}                            
\def\RANG_S{S}                            
\def\CGPARQ{Q}                            
\def\BARIT#1{{\bar {#1}}}
\def\BARH{\BARIT{H}}
\def\LATT{{\Lambda}_N}
\def\SIGMA{{\mathcal{S}_N}}
\def\LATTC{{\bar{\Lambda}_{M}}}
\def\SIGMAC{\bar{\mathcal{S}}_{M}} 

\def\JNOT{V}

\def\CUBE{C}
\def\COP{\mathbf{F}}

\def\BARH{\bar H}

\def\RELENTR{\mathcal{R}}

\def\Oo{\mathcal{O}}

\def\R{\mathbb{R}}

\def\Z{\mathbb{Z}}

\def\EXP#1{e^{#1}}

\def\EXPECT{{\mathbb{E}}}

\def\SPACE{\;\;\;}          
\def\COMMA{\,,}             
\def\PERIOD{\,.}            
\def\SEP{{\,|\,}}           

\def\VIZ#1{(\ref{#1})}      

\def\BIGO{\Oo}

\def\qed{$\Box$}
\def\PROOF{\noindent{\sc Proof:}\ }

\def\HALF{\frac{1}{2}}

\def\GRADV{\nabla V}

\newtheorem{remark}{Remark}[section]

%% file: kprt2.17.bbl
\begin{thebibliography}{10}

\bibitem{briels}
Reinier L.~C. Akkermans and W.~J. Briels.
\newblock Coarse-grained interactions in polymer melts: A variational approach.
\newblock {\em J. Chem. Phys.}, 115(13):6210--6219, 2001.

\bibitem{AKPR}
S.~Are, M.~A. Katsoulakis, P.~Plech{\'a}{\v{c}}, and L.~Rey-Bellet.
\newblock Multibody interactions in coarse-graining schemes for extended
  systems.
\newblock {\em SIAM J. Sci. Comput.}, 31(2):987--1015, 2008.

\bibitem{bco}
L.~Bertini, E.~N.~M. Cirillo, and E.~Olivieri.
\newblock Renormalization-group transformations under strong mixing conditions:
  {G}ibbsianness and convergence of renormalized interactions.
\newblock {\em J. Statist. Phys.}, 97(5-6):831--915, 1999.

\bibitem{BZ2}
A.~Bovier and M.~Zahradn{\'{\i}}k.
\newblock A simple inductive approach to the problem of convergence of cluster
  expansions of polymer models.
\newblock {\em J. Statist. Phys.}, 100(3-4):765--778, 2000.

\bibitem{presutti}
M.~Cassandro and E.~Presutti.
\newblock Phase transitions in {I}sing systems with long but finite range
  interactions.
\newblock {\em Markov Process. Related Fields}, 2(2):241--262, 1996.

\bibitem{CKV1}
A.~Chatterjee, M.~Katsoulakis, and D.~Vlachos.
\newblock Spatially adaptive lattice coarse-grained {M}onte {C}arlo simulations
  for diffusion of interacting molecules.
\newblock {\em J. Chem. Phys.}, 121(22):11420--11431, 2004.

\bibitem{CKV2}
A.~Chatterjee, M.~Katsoulakis, and D.~Vlachos.
\newblock Spatially adaptive grand canonical ensemble {M}onte {C}arlo
  simulations.
\newblock {\em Phys. Rev. E}, 71, 2005.

\bibitem{dionrev07}
A.~Chatterjee and D.G. Vlachos.
\newblock An overview of spatial microscopic and accelerated kinetic monte
  carlo methods.
\newblock {\em J. Comput-Aided Mater. Des.}, 14(2):253--308, 2007.

\bibitem{sinno}
Jianguo Dai, W.~D. Seider, and T.~Sinno.
\newblock Coarse-grained lattice kinetic {M}onte {C}arlo simulation of systems
  of strongly interacting particles.
\newblock {\em J. Chem. Phys.}, 128(19):194705, 2008.

\bibitem{ds3}
R.~L. Dobrushin and S.~B. Shlosman.
\newblock Completely analytical interactions: constructive description.
\newblock {\em J. Statist. Phys.}, 46(5-6):983--1014, 1987.

\bibitem{g97}
E.~Espanol, M.~Serrono, and Zuniga.
\newblock Coarse-grainiing of a fluid and its relation with dissipasive
  particle dynamics and smoothed particle dynamics.
\newblock {\em Int. J. Modern Phys. C}, 8(4):899--908, 1997.

\bibitem{g98}
P.~Espanol and P.~Warren.
\newblock Statistics-mechanics of dissipative particle dynamics.
\newblock {\em Europhys. Lett.}, 30(4):191--196, 1995.

\bibitem{fierro}
Francesca Fierro and Andreas Veeser.
\newblock On the a posteriori error analysis for equations of prescribed mean
  curvature.
\newblock {\em Math. Comp.}, 72(244):1611--1634, 2003.

\bibitem{doi}
H.~Fukunaga, J.~Takimoto, and M.~Doi.
\newblock A coarse-graining procedure for flexible polymer chains with bonded
  and nonbonded interactions.
\newblock {\em J. Chem. Phys.}, 116(18):8183--8190, 2002.

\bibitem{Golden}
N.~Goldenfeld.
\newblock {\em Lectures on Phase Transitions and the Renormalization Group},
  volume~85.
\newblock Addison-Wesley, New York, 1992.

\bibitem{Hadji}
G.~Hadjipanayis, editor.
\newblock {\em Magnetic Hysteresis in Novel Magnetic Materials}, volume 338 of
  {\em NATO ASI Series E}, Dordrecht, The Netherlands, 1997. Kluwer Academic
  Publishers.

\bibitem{vagelis}
V.A. Harmandaris, N.P. Adhikari, N.F.A. van~der Vegt, and K.~Kremer.
\newblock Hierarchical modeling of polystyrene: From atomistic to
  coarse-grained simulations.
\newblock {\em Macromolecules}, 39:6708--6719, 2006.

\bibitem{Kad}
L.~Kadanoff.
\newblock Scaling laws for {I}sing models near $t_c$.
\newblock {\em Physics}, 2:263, 1966.

\bibitem{KMV}
M.~A. Katsoulakis, A.~J. Majda, and D.~G. Vlachos.
\newblock Coarse-grained stochastic processes and {M}onte {C}arlo simulations
  in lattice systems.
\newblock {\em J. Comp. Phys.}, 186:250--278, 2003.

\bibitem{KPRT2}
M.~A. Katsoulakis, L.~Rey-Bellet, P.~Plech\'a\v{c}, and D.~K.Tsagkarogiannis.
\newblock Mathematical strategies in the coarse-graining of extensive systems:
  error quantification and adaptivity.
\newblock {\em J. Non Newt. Fluid Mech.}, 152:101--112, 2008.

\bibitem{KPR}
Markos~A. Katsoulakis, Petr Plech{\'a}{\v{c}}, and Luc Rey-Bellet.
\newblock Numerical and statistical methods for the coarse-graining of
  many-particle stochastic systems.
\newblock {\em J. Sci. Comput.}, 37(1):43--71, 2008.

\bibitem{KPRT}
Markos~A. Katsoulakis, Petr Plech{\'a}{\v{c}}, Luc Rey-Bellet, and Dimitrios~K.
  Tsagkarogiannis.
\newblock Coarse-graining schemes and a posteriori error estimates for
  stochastic lattice systems.
\newblock {\em M2AN Math. Model. Numer. Anal.}, 41(3):627--660, 2007.

\bibitem{KPS}
Markos~A. Katsoulakis, Petr Plech{\'a}{\v{c}}, and Alexandros Sopasakis.
\newblock Error analysis of coarse-graining for stochastic lattice dynamics.
\newblock {\em SIAM J. Numer. Anal.}, 44(6):2270--2296, 2006.

\bibitem{KT}
Markos~A. Katsoulakis and Jos{\'e} Trashorras.
\newblock Information loss in coarse-graining of stochastic particle dynamics.
\newblock {\em J. Stat. Phys.}, 122(1):115--135, 2006.

\bibitem{kremerplathe}
K.~Kremer and F.~Muller-Plathe.
\newblock Multiscale problems in polymer science: Simulation approaches.
\newblock {\em MRS Bulletin}, 26(3):205--210, 2001.

\bibitem{nochetto}
O.~Lakkis and R.~H. Nochetto.
\newblock A posteriori error analysis for the mean curvature flow of graphs.
\newblock {\em SIAM J. Numer. Anal.}, 42(5):1875--1898, 2005.

\bibitem{binder}
D.~Landau and K.~Binder.
\newblock {\em A Guide to {M}onte {C}arlo Simulations in Statistical Physics}.
\newblock Cambridge University Press, 2000.

\bibitem{g219}
A.~P. Lyubartsev, M.~Karttunen, P.~Vattulainen, and A.~Laaksonen.
\newblock On coarse-graining by the inverse monte carlo method: Dissipative
  particle dynamics simulations made to a precise tool in soft matter modeling.
\newblock {\em Soft Materials}, 1(1):121--137, 2003.

\bibitem{mp}
F.~M\"uller-Plathe.
\newblock Coarse-graining in polymer simulation: from the atomistic to the
  mesoscale and back.
\newblock {\em Chem. Phys. Chem.}, 3:754, 2002.

\bibitem{olivieri}
E.~Olivieri.
\newblock On a cluster expansion for lattice spin systems: a finite-size
  condition for the convergence.
\newblock {\em J. Statist. Phys.}, 50(5-6):1179--1200, 1988.

\bibitem{op}
E.~Olivieri and P.~Picco.
\newblock Cluster expansion for {$d$}-dimensional lattice systems and
  finite-volume factorization properties.
\newblock {\em J. Statist. Phys.}, 59(1-2):221--256, 1990.

\bibitem{PK06}
I.~Pivkin and G.~Karniadakis.
\newblock Coarse-graining limits in open and wall-bounded dissipative particle
  dynamics systems.
\newblock {\em J. Chem. Phys.}, 124:184101, 2006.

\bibitem{plass}
R.~Plass, J.A. Last, N.C. Bartelt, and G.L. Kellogg.
\newblock Self-assembled domain patterns.
\newblock {\em Nature}, 412:875, 2001.

\bibitem{clementi}
M.~Praprotnik, S.~Matysiak, L.~Delle Site, K.~Kremer, and C.~Clementi.
\newblock Adaptive resolution simulation of liquid water.
\newblock {\em J. Physics: Condensed Matter}, 19(29):292201 (10pp), 2007.

\bibitem{andelman}
M.~Seul and D.~Andelman.
\newblock Domain shapes and patterns: the phenomenology of modulated phases.
\newblock {\em Science}, 267:476--483, 1995.

\bibitem{Simon}
B.~Simon.
\newblock {\em The Statistical Mechanics of Lattice Gases, vol. {I}}.
\newblock Princeton series in Physics, 1993.

\bibitem{TT}
J.~Trashorras and D.~K. Tsagkarogiannis.
\newblock Reconstruction schemes for coarse-grained stochastic lattice systems.
\newblock {\em preprint}, 2008.
\newblock submitted.

\bibitem{tsop}
W.~Tsch\"op, K.~Kremer, O.~Hahn, J.~Batoulis, and T.~B\"urger.
\newblock Simulation of polymer melts. {II}. from coarse-grained models back to
  atomistic description.
\newblock {\em Acta Polym.}, 49:75, 1998.

\bibitem{voth}
G.A. Voth.
\newblock {\em Coarse-Graining of Condensed Phase and Biomolecular Systems}.
\newblock CRC Press, Boca Raton, FL, 2009.

\end{thebibliography}
